\documentclass[final]{svjour3PP}

\usepackage{a4,amsfonts,amsmath,amssymb,selinput,graphicx}

\usepackage{mathptmx}

\AtBeginDocument{%
  \setlength{\oddsidemargin}{\dimexpr(\paperwidth-\textwidth)/2-1in}%
  \setlength{\evensidemargin}{\oddsidemargin}%
  \setlength{\topmargin}{%
    \dimexpr(\paperheight-\textheight)/2-\headheight-\headsep-1in}%
}

\def\0{\ensuremath{\mathbf{0}}}

\def\be{\begin{equation}}
\def\bea{\begin{eqnarray*}}
\def\ee{\end{equation}}
\def\eea{\end{eqnarray*}}
\def\ba{\begin{array}}
\def\ea{\end{array}}
\def\bi{\begin{itemize}}
\def\ei{\end{itemize}}
\newtheorem{theo}{Theorem}

\newtheorem{cor}{Corollary}

\title{Stochastic subspace correction methods and fault tolerance}
\author{Michael Griebel \and Peter Oswald}
\institute{M. Griebel\at Institute for Numerical Simulation, Universit\"at Bonn, Wegelerstr. 6, 53115 Bonn, and Fraunhofer Institute for Algorithms and Scientific Computing (SCAI), Schloss Birlinghoven, 53754 Sankt Augustin\\
\\
Corresponding author, tel.: +49-228-733437, fax: +49-228-737527,\\ 
\email{griebel@ins.uni-bonn.de}
\and P. Oswald \at Institute for Numerical Simulation, Universit\"at Bonn, Wegelerstr. 6, 53115 Bonn,\\
\email{agp.oswald@gmail.com}}

\titlerunning{Stochastic Schwarz and fault tolerance}
\authorrunning{M. Griebel \and P. Oswald}

\date{}
\begin{document}
\maketitle
\begin{abstract}
We present convergence results in expectation for stochastic subspace correction schemes and their accelerated versions to solve symmetric positive-definite variational problems, and discuss their potential  for achieving fault tolerance in an unreliable compute network. We employ the standard overlapping domain decomposition algorithm for PDE discretizations to discuss the latter aspect.

\end{abstract}

\section{Introduction}\label{sec0}
With the advent of petascale compute systems in the recent years and exascale computers to arrive in the near future, there is tremendous parallel compute power available for huge parallel simulations. While this technological development surely further empowers numerical simulation as a third way of science besides theory and lab experiment, it also poses challenges: Huge systems with hundreds of thousand or even  millions of processor units will be more and more prone to failures which can corrupt the results of parallel solvers or renders them obsolete at all. It is predicted that large parallel applications may suffer from faults as frequently as once every 30 minutes on future exascale platforms \cite{SW2014}. Thus, for growing parallel computers there is the need to develop not just scalable and fast parallel algorithms but to make them fault-tolerant as well. 
Besides hard errors for which hardware mitigation techniques are under development, there is the issue of soft errors, which are either software-based or algorithm-based.
We refer to \cite{KP2013,SW2014,Tr2005} for further information and references. 

In this paper, we focus on algorithm-based fault tolerance (ABFT) and discuss how to make standard scalable and parallelizable algorithms such as domain decomposition methods in applications to partial differential equations (PDE) more fault tolerant. 
To this end, we consider stochastic subspace correction algorithms and develop a general theoretical foundation for their convergence rates  under independence assumptions for the random failure of subproblem solves.
The attractive feature of stochastic subspace correction schemes in this respect
is the fact that hard faults such as compute node failure or communication losses (as long as they are detectable) can be modeled as a random process of selecting the set of acceptable subproblem solves in each iteration step. This random process fits the independence assumption of our theory which thus allows to obtain convergence rates in the presence of faults.
We use a standard overlapping domain decomposition (DD) method for a simple two-dimensional Poisson problem as prototypical example of scalable and asymptotically optimal subspace correction methods for solving second-order elliptic PDE problems. Altogether, this shows that our convergence theory for parallel stochastic subspace correction methods indeed gives proven convergence rates also in the faulty case and results in the design of fault-tolerant methods in this setting.

The remainder of this paper is organized as follows. 
In section \ref{sec1}, we present our theoretical findings on general stochastic subspace correction methods for elliptic PDEs. For randomly chosen sets of subspaces in the iterative solution method, we prove  a bound for the convergence rate in expectation for both, a conventional one-step Schwarz-type method and its two-step Nesterov-type counterpart.
Then, in section \ref{sec11}, we consider the case of an overlapping domain decomposition method with additional coarse grid problem as a special example of our theory and discuss two types of compute systems this parallel algorithm might run on. First, we analyze a simple master-slave network were copies of data are kept on a reliable master node and the slave nodes are executing the parallel computation but are allowed to be faulty. Then, we deal with a local communication architecture which avoids global storage and global communication as much as possible and employs decentralized data storage. We give the corresponding bounds for the error reduction per iteration step in the faulty setting where the failure arrival and the failure length times are modeled by Weibull distribution functions. Moreover, we provide estimates for the parallel cost complexities of these methods.
Finally we give some concluding remarks in section \ref{sec4}. Details on a distributed implementation of the algorithm for the local communication case are discussed in the appendix.

\section{Theoretical results}\label{sec1}
Let $V$ be a separable real Hilbert space with scalar product $(\cdot,\cdot)_V$.   For a given continuous and coercive Hermitian form $a(\cdot,\cdot)$ on $V$ and a bounded linear functional $F$ on $V$, we consider the variational problem of finding the unique element $u\in V$ such that
\be\label{VP}
a(u,v)=F(v)\qquad\forall v\in V.
\ee
Equivalently, (\ref{VP}) can be formulated as quadratic minimization problem in $V$ (or as linear operator equation in the dual space $V^\ast$). In the following, we use the fact that 
$a(\cdot,\cdot)$ defines a spectrally equivalent scalar product on $V$, equip $V$ with it, and write $\|v\|=a(v,v)^{1/2}$. In our model application, this $V$ will be a large FE subspace of $H^1_0(\Omega)$, the discretization space of a diffusion problem
\be\label{MP}
-\nabla_x\cdot(a(x)\nabla_x u(x))=f(x),\qquad u\in H^1_0(\Omega), \quad 0\le c_0 \le a(x) \le C_0,\quad x\in \Omega,
\ee
where $\Omega$ is a (nice) domain in $\mathbb{R}^d$ for moderate dimension $d$. The associated variational 
problem reads
$$
a(u,v):=\int_\Omega a(x)\nabla_x u(x)\cdot \nabla_xv(x)\, dx = F(v):= \int_\Omega f(x)v(x)\,dx.
$$

We consider Schwarz iterative methods (also called subspace correction methods) for solving (\ref{VP}). The underlying {\it space splitting} is given by
a collection $\{V_i\}_{i=0,1,...,n}$ of $n+1$ separable real Hilbert spaces, each equipped with a spectrally equivalent scalar product 
$a_i(\cdot,\cdot)$ and norm $a_i(v_i,v_i)^{1/2}$, and bounded linear operators $R_i:\,V_i\to V$ such that 
$$
\sum_{i=0}^n R_iV_i = V,
$$
and $\mathrm{Ker}(R_i)=\{0\}$ for all $i=0,1,\ldots,n$ (otherwise replace $V_i$ by $V_i \ominus_i \mathrm{Ker}(R_i)$). 
This induces another set of bounded linear operators 
$T_i=R_i^\ast :\,V\to V_i$, $i=0,1,\ldots,n$, defined by the solution of auxiliary variational problems in $V_i$:
\be\label{VPo}
a_i(T_i v,v_i)=a(v,R_i v_i)\qquad \forall\;v_i\in V_i.
\ee
In the $m$-th step of a Schwarz iterative method for solving (\ref{VP}), a certain finite set $I_m\subset \{0,1,\ldots,n\}$ is chosen (deterministically, randomly, or in a greedy fashion), for each $i\in I_m$ the corresponding auxiliary problem (\ref{VPo})
is solved with $v=e^{(m)}_u=u-u^{(m)}$, and an update of the form
\be\label{Rec0}
u^{(m+1)}=u^{(m)} + \sum_{i\in I_m} \xi_{m,i}R_i T_i e^{(m)}_u,\qquad m=0,1,\ldots,
\ee
is performed.  At start, we set w.l.o.g. $u^{(0)}=0$. Since
$$
a(e^{(m)}_u,R_i v_i)=F(R_i v_i) - a(u^{(m)},R_i v_i), \qquad v_i\in V_i,
$$
the update step (\ref{Rec0}) can be performed without knowing $u$.
The relaxation parameters $\xi_{m,i}$, $i=0,1,\ldots,n$,
can be chosen depending on $u^{(m)}$ (then a nonlinear iteration results), or independent of
$u^{(m)}$  (then we have a generally non-stationary but linear iterative scheme). The iteration (\ref{Rec0}) subsumes different standard algorithms such as the multiplicative (or sequential) Schwarz method
where in each step a single subproblem (\ref{VPo}) is solved ($|I_m| =1$), the additive (or parallel) Schwarz method where all $n$ subproblems are solved simultaneously ($I_m=\{0,1,\ldots,n\}$), and intermediate block-iterative schemes ($1<|I_m| <n+1$). Here and in the following, $|I|$ denotes
the cardinality of a finite index set $I$. The recursion (\ref{Rec0}) basically represents a one-step iterative
method (only the current iterate $u^{(m)}$ needs to be available for the update step). Below we will also consider an accelerated scheme by introducing a two-step Schwarz method in vector form inspired by \cite{Ne2012}. We refer to \cite{Os1994,Xu1992} for early work on Schwarz iterative methods. 

In this paper, we focus on stochastic versions of the one-step method (\ref{Rec0}) and its two-step counterpart, where the sets $I_m$ are chosen randomly.
To be precise, we assume that
\bi
\item[{\bf A}] $I_m$ is a uniformly at random
chosen subset of size $p_m$ in $\{0,1,\ldots,n\}$, i.e., $|I_m|=p_m$ and  $\mathbb{P}(i\in I_m)=\mathbb{P}(i'\in I_m)$ for all
$i,i'\in \{0,1,\ldots,n\}$. 
\item[{\bf B}] The choice of $I_m$ is independent for different $m$.
\ei
Below, we will consider expectations of squared error norms for iterations with any fixed but arbitrary sequence $\{p_m\}$.
Here, the restriction to uniform index sampling in {\bf A} is not essential, see the remarks after Theorem \ref{theo1}. What is important
for the proof technique is the independence assumption {\bf B}.

For such a random choice of $I_m$, a convergence estimate for the expectation of the squared error in terms of the stability constants
of the space splitting has already been announced without proof in \cite{GrOs2012} (see Theorem 3 in \cite{OsZh2015} for the argument). 
We formulate it in a slightly modified setting including weights for convenience. Let $0<\lambda_{\min}\le \lambda_{\max} < \infty$ and positive
weights $\omega=\{\omega_i>0\}_{i=0,1,\ldots,n}$ be such that it holds
\be\label{NE}
\lambda_{\min} \||v\||_\omega^2 \le a(v,v) \le \lambda_{\max}\||v\||_\omega^2,\qquad v\in V,
\ee
where
$$
\||v\||_\omega^2:=\inf_{\ba{c}v_i\in V_i,\,i=0,1,\ldots,n \\ \;v=\sum_{i=0}^n \omega_i R_i v_i\ea} \sum_{i=0}^n \omega_i a_i(v_i,v_i).
$$
The norm equivalence (\ref{NE}) can also be written in terms of properties of the operators $T_i$, $R_i$, and the additive Schwarz operator
$$
P :=\sum_{i=0}^n\omega_i R_iT_i \,: V\to V
$$
associated with the space splitting as follows: It holds
\be\label{NE1}
\lambda_{\min} a(v,v)  \le a(Pv,v)=\sum_{i=0}^n \omega_i a_i(T_iv,T_iv), \qquad v\in V,
\ee
and
\be\label{NE2}
 \|\sum_{i=0}^n \omega_i R_iv_i\|^2  \le 
 \lambda_{\max}\sum_{i=0}^n \omega_i a_i(v_i,v_i),\qquad v_i\in V_i.
\ee
It is well known that the stability condition (\ref{NE})  of the space splitting implies that $P$
is positive-definite with respect to $a(\cdot,\cdot)$, and satisfies
\be\label{NE3}
\||v\||_\omega^2 = a(P^{-1}v,v),\qquad v\in V,\qquad\quad\lambda_{\min}  \mathrm{Id} \le P  \le  \lambda_{\max} \mathrm{Id}.
\ee
We now can state the following convergence result:

\begin{theo}\label{theo1}
Let the relaxation parameters in (\ref{Rec0}) be given by $\xi_{m,i}:=\xi \omega_i$, $i=0,1,\ldots,n$,
where $0<\xi < 2/ \lambda_{\max}$. Furthermore, let the random sets $I_m$ of size $p_m$ be selected in agreement with {\bf A}. Then in each step the algorithm (\ref{Rec0}) reduces the error in expectation according to  
\be\label{EC}
\mathbb{E}(\|e_u^{(m+1)}\|^2\,|\,u^{(m)}) \le \left(1-\frac{ \lambda_{\max}\xi(2-  \lambda_{\max}\xi) p_m}{\kappa (n+1)}\right) \|e_u^{(m)}\|^2,\qquad m=0,1,\ldots,
\ee
where $\kappa:=\lambda_{\max}/\lambda_{\min}$ is the condition number of the underlying space splitting.\\
If in addition {\bf B} holds then the algorithm (\ref{Rec0}) converges in expectation for any $u\in V$ and
\be\label{ECC}
\mathbb{E}(\|e_u^{(m)}\|^2) \le \prod_{s=0}^{m-1}\left(1-\frac{ \lambda_{\max}\xi(2-  \lambda_{\max}\xi) p_s}{\kappa (n+1)}\right) \|u\|^2,\qquad m=1,2,\ldots.
\ee
\end{theo}

{\bf Proof}. Obviously, (\ref{ECC}) directly follows from (\ref{EC}) by the independence assumption {\bf B}. Thus, it suffices to consider a single step of
the iteration, and to compute the expectation of the squared error
$\|e_u^{(m+1)}\|^2$ conditioned on an arbitrarily fixed $u^{(m)}$ (to keep the formulas short, in the proof we simply write $\mathbb{E}(\cdot)$ for conditional expectations $\mathbb{E}(\cdot\,|\,u^{(m)})$).
By (\ref{Rec0}) we have
\bea
\mathbb{E}(\|e_u^{(m+1)}\|^2) &=& \mathbb{E}(\|e_u^{(m)}\|^2-2\xi
a(e_u^{(m)},\sum_{i\in I_m} \omega_i R_iT_ie_u^{(m)}) +\xi^2
\|\sum_{i\in I_m} \omega_i R_iT_ie_u^{(m)}\|^2)\\
&=& \|e_u^{(m)}\|^2-2\xi
\mathbb{E}(a(e_u^{(m)},\sum_{i\in I_m} \omega_i R_iT_ie_u^{(m)})) +\xi^2
\mathbb{E}(\|\sum_{i\in I_m} \omega_i R_iT_ie_u^{(m)}\|^2).
\eea
For the second term, we have
$$
\mathbb{E}(a(e_u^{(m)},\sum_{i\in I_m} \omega_i R_iT_ie_u^{(m)}))=
a(e_u^{(m)},\mathbb{E}(\sum_{i\in I_m}\omega_i R_iT_ie_u^{(m)}))=
\frac{p_m}{n+1} a(e_u^{(m)},Pe_u^{(m)}),
$$
since by assumption {\bf A} and the definition of $P$ 
$$
\mathbb{E}(\sum_{i\in I_m} \omega_i R_iT_ie_u^{(m)})=
\frac{p_m}{n+1} \sum_{i=0}^n\omega_i R_iT_ie_u^{(m)}=
\frac{p_m}{n+1} Pe_u^{(m)}.
$$
This equality and (\ref{NE2}) imply for the last term that  
\bea
\mathbb{E}(\|\sum_{i\in I_m} \omega_i R_iT_ie_u^{(m)}\|^2)&\le&
\lambda_{\max} \mathbb{E}(\sum_{i\in I_m} \omega_i a_i(T_ie_u^{(m)},T_ie_u^{(m)}))
\\
&=& \lambda_{\max} a(e_u^{(m)},\mathbb{E}(\sum_{i\in I_m} \omega_i R_iT_ie_u^{(m)}))=\frac{\lambda_{\max}p_m}{n+1} a(e_u^{(m)},Pe_u^{(m)}).
\eea
Substitution gives
$$
\mathbb{E}(\|e_u^{(m+1)}\|^2)\le \|e_u^{(m)}\|^2 -\xi\frac{p_m}{n+1}(2-\lambda_{\max}\xi)a(e_u^{(m)},Pe_u^{(m)}),
$$
and the lower spectral bound in (\ref{NE3}) finally yields
$$
\mathbb{E}(\|e_u^{(m+1)}\|^2)\le (1-\lambda_{\min}\xi\frac{p_m}{n+1}(2-\lambda_{\max}\xi))\|e_u^{(m)}\|^2 = (1-\frac{\lambda_{\max}\xi(2-\lambda_{\max}\xi)p_m}{\kappa (n+1)})\|e_u^{(m)}\|^2.
$$
This proves the statement of Theorem \ref{theo1}.\hfill $\Box$

\medskip
The recent paper \cite{GrOs2017} contains similar results for infinite-dimensional $V$ and
countable splittings albeit with weaker non-geometric convergence rates in expectation under certain smoothness assumptions
on $u$.

An application of the iteration (\ref{Rec0}) with theoretical guarantees according to
Theorem \ref{theo1} requires
knowledge of suitable weights $\omega_i$, and  an upper bound $\bar{\lambda}$ for the stability constant $\lambda_{\max}$ in order to choose the value of $\xi$, whereas information about $p_m$, the size of $I_m$, is not crucial, see the discussion below. Numerical experiments for model problems with different values  $\xi\in  (0,2/\lambda_{\max})$ suggest that the iteration count is sensitive to the choice of $\xi$ and that overrelaxation gives often better results. An alternative, especially in cases when no reliable information on $\lambda_{\max}$ is available, is to choose
$\xi=\xi_m$ depending on $u^{(m)}$ and $d^{(m)}=\sum_{i\in I_m} \omega_i R_iT_ie^{(m)}$ by the steepest decent rule,
i.e., to minimize $\|e_u^{(m+1)}\|$ for given $u^{(m)}$ and $d^{(m)}$ by setting
\be\label{XiOpt}
\xi_m:= \frac{a(e^{(m)},d^{(m)})}{a(d^{(m)},d^{(m)})}.
\ee
But note that, in an implementation of (\ref{Rec0}) with the steepest descent rule for the compute networks considered in the next section,  the global scalar products needed for (\ref{XiOpt}) represent a bottleneck similar to the solution of the coarse subproblem in e.g. a DD method and the global error computation. On the positive side, due to the minimization property 
of the steepest descent rule and the method of proof for Theorem \ref{theo1}, the expectation of the squared error 
for (\ref{Rec0}) with the non-stationary steepest decent rule $\xi=\xi_m$ must satisfy the same upper bound as
the best bound with fixed $\xi$.

The weights $\omega_i$ can be considered as scaling parameters that can be used to
improve the stability constants $\lambda_{\max}$, $\lambda_{\min}$, and thus the condition number $\kappa$ of the space splitting. 
Note that improving the estimate (\ref{EC}) in Theorem \ref{theo1} by minimizing $\kappa$ via the choice of optimal weights $\omega_i$ is obviously equivalent to optimizing the set of relaxation parameters $\{\xi_{m,i}=\xi \omega_i\}$. Similar improvements can also be achieved, at least approximately,
by adapting the probability distribution for choosing the sets $I_m$. Indeed, instead of {\bf A}, assume that $I_m$ is a randomly chosen index set of size $p_m\le n+1$ such that $\mathbb{P}(i\in I_m)=p_mq_i>0$, where $\{q_i\}_{i=0,1,\ldots,n}$ is an arbitrary discrete probability distribution with support $\{0,1,\ldots,n\}$ (assumption {\bf A} corresponds to the case of a uniform distribution with $q_i=1/(n+1)$). Such a condition can, in general, be achieved only approximately (or exactly if we allow for repetitions in $I_m$). This changes a few lines in the proof of Theorem \ref{theo1}, namely, we have
$$
\mathbb{E}(a(e_u^{(m)},\sum_{i\in I_m} \omega_i R_iT_ie_u^{(m)}))
= p_m a(e_u^{(m)},\sum_{i=0}^n q_i\omega_i R_iT_ie_u^{(m)}) =
\frac{p_m}{n+1} a(e_u^{(m)}, \tilde{P} e_u^{(m)}),
$$
and, similarly,
$$
\mathbb{E}(\|\sum_{i\in I_m} \omega_i R_iT_ie_u^{(m)}\|^2)\le \lambda_{\max}p_m a(e_u^{(m)},\sum_{i=0}^n q_i\omega_i R_iT_ie_u^{(m)})=
\lambda_{\max}p_m a(e_u^{(m)},\tilde{P}e_u^{(m)}),
$$
where $\tilde{P}$ is the additive Schwarz operator for the space splitting with the same $V_i$ and $R_i$ but with a different set of weights $\tilde{\omega}_i:=q_i\omega_i$.
With these changes and by choosing $\xi=\lambda_{\max}^{-1}$, the reduction factor for the expectation of the squared error in (\ref{EC}) becomes 
$$
1-\frac{ \tilde{\lambda}_{\min}p_m}{\lambda_{\max}},
$$
where $\tilde{\lambda}_{\min}$ is the lower spectral bound associated with $\tilde{P}$. Thus, varying $q_i>0$ under the normalization
condition
$$
\sum_{i=0}^n q_i =1
$$
may also result in improved bounds.

The discussion on applications of our convergence estimates to fault tolerance below will focus on the situation when
$$
1 << p_m\le p\le n+1,
$$
where $p$ stands for the number of processors available for subproblem solves in a compute network, and $p_m$ is a sequence of random integers standing for the number of correctly working processors. In such a case, the average reduction of the expectation of the squared error per iteration step is approximately given by
$$
\left(\prod_{s=0}^{m-1} (1-\frac{p_s}{(n+1)\kappa})\right)^{1/m}
\approx 1-\frac{\sum_{s=0}^{m-1} p_s}{(m(n+1)\kappa} \approx
1-\frac{r_p}{\kappa},\qquad r_p:=\mathbb{E}(p_m)/(n+1),
$$
if we set $\xi=1/{\lambda}_{\max}$ and take sufficiently large $m$.
The number $r_p$ can be interpreted as the average rate of subproblem solves
per iteration step (\ref{Rec0}), and our estimate of the average error reduction per step suggests that the convergence speed of the recursion (\ref{Rec0}) is tied to it in a linear fashion which is as good as one can hope for.

The deterioration of the convergence rate with the condition  number $\kappa$ of the associated weighted
space splitting is typical for one-step iterations such as (\ref{Rec0}). The convergence rate can be improved to a dependence on only $\kappa^{1/2}$ rather than on $\kappa$ by using multi-step strategies.
This has recently attracted attention in the optimization literature, see e.g. \cite{FerRic2016,LeeSid2013,Ne2012,RiTa2017}.
Following \cite{LiuWr2016,Ne2012} we will consider the subsequent accelerated Schwarz method written in vector form,
as a one-step iteration for two sequences $u^{(m)}$ and $v^{(m)}$:
With $u^{(0)}=v^{(0)}=0$ at start, for $m=0,1,\ldots$ execute
\be\label{u}
u^{(m+1)} = w^{(m)} + \xi_m \sum_{i\in I_m} \omega_i R_iT_ie_w^{(m)},\qquad w^{(m)} = \alpha_m v^{(m)} + (1-\alpha_m) u^{(m)},
\ee
\be\label{v}
v^{(m+1)} = \beta_m v^{(m)}+(1-\beta_m)w^{(m)} + \eta_m \sum_{i\in I_m} \omega_i R_iT_ie_w^{(m)},\qquad\qquad\quad
\ee
with parameter sequences $\alpha_m,\beta_m,\xi_m,\eta_m$ determined below. The random index set $I_m$ is chosen according to rule {\bf A}. Notation for errors is as usual:
$$
e_u^{(m)}:=u-u^{(m)},\quad e_v^{(m)}:=u-v^{(m)},
\quad e_w^{(m)}:=u-w^{(m)}=\alpha_m e_v^{(m)}+(1-\alpha_m) e_u^{(m)}.
$$
Note that, as for the iteration (\ref{Rec0}), in each step $p_m$ subproblems have to be solved
but storage and update work increase. Remedies are available, see \cite{FerRic2016,LiuWr2016} for discussions on
implementation issues. We have the following convergence result:

\begin{theo}\label{theo2}
Assume that we possess upper and lower bounds 
$$
0<\underline{\lambda} \le \lambda_{\min}\le \lambda_{\max}\le \bar{\lambda}
$$ 
for the stability constants of the space decomposition. Then, with 
$\xi_m=\bar{\lambda}^{-1}$, $\eta_m=(\bar{\lambda}\underline{\lambda})^{-1/2}$ 
and sequences $\alpha_m,\beta_m\in (0,1)$ defined in (\ref{AM}), (\ref{BM}) below and under the assumption {\bf A},  the vector iteration (\ref{u}-\ref{v}) admits the recursive estimate
\be\label{ECA}
\mathbb{E}(\|e_u^{(m+1)}\|^2 + \underline{\lambda}\||e_v^{(m+1)}\||_\omega^2\,|\,u^{(m)},v^{(m)})\le \left(1 - \frac{p_m}{(n+1)\sqrt{\bar{\kappa}}}\right)(\|e_u^{(m)}\|^2 + \underline{\lambda}\||e_v^{(m)}\||_\omega^2) ,\qquad 
\ee
$m=0,1,\ldots$. Here, $\bar{\kappa}=\bar{\lambda}/\underline{\lambda}$ is an upper bound for the condition $\kappa$ of the space splitting.\\
If in addition {\bf B} holds then the algorithm (\ref{u}-\ref{v}) converges in expectation for any $u\in V$, and
\be\label{ECAU}
\mathbb{E}(\|e_u^{(m)}\|^2) \le 2\prod_{s=0}^{m-1}\left(1 - \frac{p_s}{(n+1)\sqrt{\bar{\kappa}}}\right)\|u\|^2,\qquad m=1,2,\ldots.
\ee
\end{theo}

{\bf Proof}. The proof is a slightly simplified adaption of the proofs in \cite{LiuWr2016,Ne2012} which takes into account that
we allow for redundant space splittings and variable block sizes $p_m$ in our iteration method. In the proof of (\ref{ECA}), we again shorten notation, and write $\mathbb{E}(\cdot)$ instead of $\mathbb{E}(\cdot\,|\,u^{(m)},v^{(m)})$. 

From (\ref{u}) we compute with $\xi_m=\bar{\lambda}^{-1}$
$$
\|e_u^{(m+1)}\|^2 =\|e_w^{(m)}\|^2 - 2\bar{\lambda}^{-1}a(e_w^{(m)},\sum_{i\in I_m} \omega_i R_iT_ie_w^{(m)})+\bar{\lambda}^{-2}\|\sum_{i\in I_m} \omega_i R_iT_ie_w^{(m)}\|^2.
$$
When we take the conditional expectation with respect to $I_m$ and recall the definition and the properties of the additive Schwarz operator $P$, we get due to {\bf A}
$$
\mathbb{E}(a(e_w^{(m)},\sum_{i\in I_m} \omega_i R_iT_ie_w^{(m)})) = \frac{p_m}{n+1} a(e_w^{(m)},\sum_{i=0}^n \omega_i R_iT_ie_w^{(m)}) = 
\frac{p_m}{n+1} a(Pe_w^{(m)},e_w^{(m)}),
$$
and by (\ref{NE2}) we obtain
\bea
\mathbb{E}(\|\sum_{i\in I_m} \omega_i R_iT_ie_w^{(m)}\|^2)&\le& \lambda_{\max}\mathbb{E}(\sum_{i\in I_m} \omega_i a_i(T_ie_w^{(m)},T_ie_w^{(m)})) = 
\lambda_{\max}a(\mathbb{E}(\sum_{i\in I_m} \omega_i R_iT_ie_w^{(m)}),e_w^{(m)})\\
&=&
\frac{\lambda_{\max}p_m}{n+1} a(Pe_w^{(m)},e_w^{(m)}) \le \frac{\bar{\lambda}p_m}{n+1} a(Pe_w^{(m)},e_w^{(m)}).
\eea
This gives
$$
\mathbb{E}(\|e_u^{(m+1)}\|^2)\le \|e_w^{(m)}\|^2 - \frac{p_m}{(n+1)\bar{\lambda}}a(Pe_w^{(m)},e_w^{(m)}),
$$
or after rearrangement
\be\label{EU}
\frac{p_m}{n+1}a(Pe_w^{(m)},e_w^{(m)})\le \bar{\lambda}(\|e_w^{(m)}\|^2 - \mathbb{E}(\|e_u^{(m+1)}\|^2)).
\ee
For the estimation of the error term $\||e_v^{(m+1)}\||_\omega^2$ recall (\ref{NE3}). We can now write 
\bea 
\||e_v^{(m+1)}\||_\omega^2&=&\||\beta_m e_v^{(m)} + (1-\beta_m)e_w^{(m)}\||_\omega^2\\
&&\; -2\eta_m a(P^{-1}(\beta_m e_v^{(m)} + (1-\beta_m)e_w^{(m)}),
\sum_{i\in I_m} \omega_i R_iT_ie_w^{(m)}) + \eta_m^2 \||\sum_{i\in I_m} \omega_i R_iT_ie_w^{(m)}\||_\omega^2.
\eea 
After taking the expectation with respect to $I_m$, each of the three terms in the right-hand side will be estimated separately. For the first term (which does not depend on $I_m$), we have by convexity of the norm
$$
\||\beta_m e_v^{(m)} + (1-\beta_m)e_w^{(m)}\||_\omega^2\le \beta_m \||e_v^{(m)}\||_\omega^2 + (1-\beta_m)\||e_w^{(m)}\||_\omega^2.
$$
Using the lower bound in (\ref{NE}) (or, equivalently, the upper bound for the spectrum of $P^{-1}$ in (\ref{NE3})), this yields
\be\label{EV1}
A_1:=\||\beta_m e_v^{(m)} + (1-\beta_m)e_w^{(m)}\||_\omega^2\le  \beta_m \||e_v^{(m)}\||_\omega^2 +
(1-\beta_m)\underline{\lambda}^{-1}\|e_w^{(m)}\|^2.
\ee
For the second term, acting as before we get
\bea
&&A_2:=\mathbb{E}(a(P^{-1}(\beta_m e_v^{(m)} + (1-\beta_m)e_w^{(m)}),\sum_{i\in I_m} \omega_i R_iT_ie_w^{(m)}))\\
&& \qquad\qquad =\frac{p_m}{n+1} a(P^{-1}(\beta_m e_v^{(m)} + (1-\beta_m)e_w^{(m)}),\sum_{i=0}^n\omega_i R_iT_ie_w^{(m)})\\
&& \qquad\qquad =\frac{p_m}{n+1} a(P^{-1}(\beta_m e_v^{(m)} + (1-\beta_m)e_w^{(m)}),Pe_w^{(m)}))\\
&& \qquad\qquad =\frac{p_m}{n+1} (\|e_w^{(m)}\|^2 + \beta_m a(e_v^{(m)}- e_w^{(m)},e_w^{(m)}).
\eea
Here, to eliminate $e_v^{(m)}$ from the last term, we use the definition of $w^{(m)}$ in (\ref{u}): Since 
$$
e_v^{(m)}-e_w^{(m)}=\frac{1-\alpha_m}{\alpha_m}(e_w^{(m)}-e_u^{(m)}),
$$
we obtain
\bea
a(e_v^{(m)}- e_w^{(m)},e_w^{(m)})&=&\frac{1-\alpha_m}{\alpha_m}a(e_w^{(m)} - e_u^{(m)},e_w^{(m)})\\&=&\frac{1-\alpha_m}{2\alpha_m}(\|e_w^{(m)}-e_u^{(m)}\|^2+\|e_w^{(m)}\|^2 -\|e_u^{(m)}\|^2)\\
&\ge & \frac{1-\alpha_m}{2\alpha_m}(\|e_w^{(m)}\|^2 - \|e_u^{(m)}\|^2).
\eea
Thus,
\be\label{EV2}
A_2\ge \frac{p_m}{n+1}((1 + \frac{\beta_m(1-\alpha_m)}{2\alpha_m})\|e_w^{(m)}\|^2 - \frac{\beta_m(1-\alpha_m)}{2\alpha_m}\|e_u^{(m)}\|^2).
\ee
Finally, for the last term we first use the definition of the $\||\cdot\||_\omega$ norm which gives
$$
\||\sum_{i\in I_m} \omega_i R_iT_ie_w^{(m)}\||_\omega^2\le \sum_{i\in I_m} \omega_i a_i(T_ie_w^{(m)},T_ie_w^{(m)})=a(\sum_{i\in I_m} \omega_i R_i T_i e_w^{(m)},e_w^{(m)},
$$
and then proceed as above to arrive at
$$
A_3:=\mathbb{E}(\||\sum_{i\in I_m} \omega_i R_iT_ie_w^{(m)}\||_\omega^2)\le \frac{p_m}{n+1} a(Pe_w^{(m)},e_w^{(m)}).
$$
It remains to use (\ref{EU}) which gives
\be\label{EV3}
A_3
\le \bar{\lambda}(\|e_w^{(m)}\|^2 - \mathbb{E}(\|e_u^{(m+1)}\|^2)).
\ee
Substitution of (\ref{EV1}-\ref{EV3}) and collecting all multiples of $\|e_w^{(m)}\|^2$ into one expression gives
\bea
\mathbb{E}(\||e_v^{(m+1)}\||_\omega^2)&=& A_1 - 2\xi_m A_2 + \xi_m^2 A_3\\
&\le& \beta_m \||e_v^{(m)}\||_\omega^2 + \frac{p_m}{n+1}\frac{\eta_m\beta_m(1-\alpha_m)}{\alpha_m}\|e_u^{(m)}\|^2 - \bar{\lambda}\eta_m^2\mathbb{E}(\|e_u^{(m+1)}\|^2\\
&&\quad +((1-\beta_m)\underline{\lambda}^{-1} - 2\eta_m\frac{p_m}{n+1}(1 + \frac{\beta_m(1-\alpha_m)}{2\alpha_m}) + \bar{\lambda}\eta_m^2)\|e_w^{(m)}\|^2.
\eea
Under the assumption that we will be able to choose the coefficient in front of $\|e_w^{(m)}\|^2$ equal to zero (or negative), i.e., if
\be\label{C1}
(1-\beta_m)\underline{\lambda}^{-1} - 2\eta_m\frac{p_m}{n+1}(1 + \frac{\beta_m(1-\alpha_m)}{2\alpha_m}) + \bar{\lambda}\eta_m^2\le 0,
\ee
this turns into the inequality
\be\label{ECA1}
\mathbb{E}(\||e_v^{(m+1)}\||_\omega^2) + \bar{\lambda}\eta_m^2\mathbb{E}(\|e_u^{(m+1)}\|^2)\le \beta_m(\||e_v^{(m)}\||_\omega^2 + \frac{p_m\eta_m(1-\alpha_m)}{(n+1)\alpha_m}\|e_u^{(m)}\|^2).
\ee
To eventually arrive at (\ref{ECA}) the remaining steps are to choose $\alpha_m,\beta_m\in (0,1)$ and $\eta_m>0$ in such a way that in addition to (\ref{C1}) 
the coefficients in front of $\mathbb{E}(\|e_u^{(m+1)}\|^2)$ and $\|e_u^{(m)}\|^2$ in (\ref{ECA1}) coincide, i.e.,
\be\label{C2}
\bar{\lambda}\eta_m^2=\frac{p_m\eta_m(1-\alpha_m)}{(n+1)\alpha_m},
\ee
and that $\beta_m$ is as small as possible. This is done as follows: The parameter $\alpha_m$ can always be determined such that (\ref{C2}) holds
while (\ref{C2}) can be used to eliminate $\alpha_m$ from (\ref{C1}). Indeed, from (\ref{C2}) we have
\be\label{AM}
\frac{1-\alpha_m}{\alpha_m}=\frac{(n+1)\bar{\lambda}}{p_m}\eta_m,
\ee
and (\ref{C1}) turns into
\bea
&&(1-\beta_m)\underline{\lambda}^{-1} - \eta_m(2\frac{p_m}{n+1}+\beta_m\bar{\lambda}\eta_m) + \bar{\lambda}\eta_m^2)\\
&&\qquad =(1-\beta_m)(\frac1{\underline{\lambda}} - \frac{2p_m}{(n+1)(1-\beta_m)}\eta_m + \bar{\lambda}\eta_m^2)\le 0.
\eea
This inequality has positive solutions $\eta_m$ iff the polynomial
$$
p(t)= \frac1{\underline{\lambda}} - \frac{2p_m}{(n+1)(1-\beta_m)}t + \bar{\lambda}t^2
$$
has a positive real root, i.e., if its discriminant is non-negative. This implies the condition
$$
\frac{p_m^2}{(n+1)^2(1-\beta_m)^2} - \frac{\bar{\lambda}}{\underline{\lambda}}=\frac{p_m^2}{(n+1)^2(1-\beta_m)^2} - \bar{\kappa}\ge 0.
$$
The smallest possible $\beta_m$ and associated $\eta_m>0$ for which this inequality may hold are thus given by
\be\label{BM}
\beta_m = 1	-\frac{p_m}{(n+1)\bar{\kappa}^{1/2}},\qquad \eta_m = (\bar{\lambda}\underline{\lambda})^{-1/2}.
\ee
Now, the coefficient in front of $\mathbb{E}(\|e_u^{(m+1)}\|^2)$ and $\|e_u^{(m)}\|^2$ in (\ref{ECA1})
equals $\bar{\lambda}\eta_m^2=\underline{\lambda}^{-1}$, and is independent of $m$. Multiplying in (\ref{ECA1}) by 
$\underline{\lambda}$ gives (\ref{ECA}). 

Due to our assumptions and (\ref{NE}) we have
$$
\underline{\lambda}\||e_v^{(m)}\||_\omega^2\le \lambda_{\min}\||e_v^{(m)}\||_\omega^2\le \|e_v^{(m)}\|^2,
$$
and $e_v^{(0)}=e_u^{(0)}=u$. Thus, together with the independence assumption {\bf B}, the per step estimate (\ref{ECA}) implies the convergence in expectation for arbitrary $u\in V$ and the bound (\ref{ECAU}) for the (unconditional) expectation
of the squared error $\|e_u^{(m)}\|^2$.
This completes the proof of Theorem \ref{theo2}.\hfill $\Box$

\medskip
That, in contrast to the recursion (\ref{Rec0}), the coefficient $\alpha_m$ in the recursion formula (\ref{u}) for the accelerated scheme depends on the size $p_m$ of the random index set $I_m$ is not a problem as long as $I_m$ is known before the subproblems needed for the update steps (\ref{u}) and (\ref{v}) are solved. However,
in the applications discussed in the next section this is not the case: The appropriate set $I_m$ is known only after
the required subproblem solves are executed. Consequently, the coefficient $\alpha_m$ which is needed to compute the residuals for the
subproblems may not have been chosen properly. Remedies are to work with safe lower bounds $\underline{p}_m$ for $p_m$ (thus giving up some accuracy in the upper bound (\ref{ECAU}) by replacing $p_m$ by $\underline{p}_m$), or to perform twice as many subproblem solves, namely to compute 
$R_iT_ie_u^{(m)}$ and $R_iT_ie_v^{(m)}$ separately for each required $i$, and to perform the linear combinations only after the exact $I_m$, and thus
$\alpha_m$, is known. 

To summarize, the above convergence theory covers a stochastic version of Schwarz iterative methods based on generic space splittings, where
in each step a random subset of subproblem solves is used. On the one hand, this theory shows that randomized Schwarz iterative
methods are competitive with their deterministic counterparts. On the other hand, there are situations where randomness in the subproblem selection is naturally occurring, and not a matter of choice in the numerical method. An important example is given by
algorithm-based methods for achieving fault tolerance in large-scale distributed and parallel computing applications. This will be dealt with in more detail in the remainder of this paper using overlapping domain decomposition PDE solvers as
an example.

\section{Potential for achieving fault tolerance}\label{sec11}
The occurrence of faults in the execution of large-scale computational tasks and their mitigation has become an issue in recent years due to cloud computing applications and the exascale HPC development. The common assumption is that hard and soft faults may occur more often in future computing applications, that they are probabilistic in nature, and that strategies to counteract them will become increasingly important.
A wide range of approaches and proposals for achieving fault tolerance with little overhead (often based on hypotheses on the fault model  applicable to future exascale computer architectures that are hard to validate at present) are currently under consideration. They are commonly categorized as hardware-based (HBFT), software-based (SBFT), and algorithm-based (ABFT). We refer to \cite{KP2013,SW2014,Tr2005} for more information and references. 

We now concentrate on the ABFT aspect, and discuss the potential consequences of  the theoretical results from section \ref{sec1} for making standard scalable and parallelizable algorithms such as domain decomposition methods 
in PDE applications more fault tolerant. The attractive feature of stochastic subspace correction schemes in this respect
is the fact that hard faults such as compute node failure or communication losses (as long as they are detectable) can be modeled as a random process of selecting the index set $I_m$ of acceptable subproblem solves in each iteration step (\ref{Rec0}). This random process often fits the independence assumption {\bf B} that is crucial in order to obtain the convergence rates in Theorem \ref{theo1}, similarly for the accelerated version (\ref{u}-\ref{v}) and Theorem \ref{theo2}. 

We use a standard overlapping domain decomposition (DD) method for the model PDE problem (\ref{MP}) with $a(x)=1$  in the following discussion since it represents one of the prototypical examples of scalable and asymptotically optimal subspace correction methods for solving second-order elliptic PDE problems. The parallelization of overlapping DD solvers is (up to the solution of the so-called coarse problem, see below) straightforward, even though most implementations are based on non-overlapping DD schemes which provide better data locality at the expense of sometimes asymptotically non-optimal preconditioning behavior, compare \cite{GZu,KG1987,Sm1993}.  Many statements we make carry over, with minor modifications, to non-overlapping DD methods and multigrid schemes for which DD-type implementations are used for their parallel execution. More details on the DD method for
(\ref{MP}) will be given in the next subsection \ref{sec12}. 

Whether fault tolerance can be achieved without significant cpu-time penalty depends very much on the compute architecture, and in particular on the relative speed difference between compute and communication steps
in distributed or parallel computer networks. To this end, we will below discuss
different scenarios to illustrate the application of our theoretical results. 

In subsection
\ref{sec13} we consider a master-slave architecture (think of outsourcing of compute effort to the cloud controlled by a reliable server with large and fast memory access) neglecting all overhead due to communication between slave nodes and master node. Under the assumption that each slave can execute one subproblem per cycle but can fail 
to return correct results with a certain probability, we present a random assignment scheme for which convergence in expectation can be guaranteed by our results of section \ref{sec1} independently of the fault process. Numerical experiments for
a generic domain decomposition scheme for the Poisson problem in two dimensions demonstrate the robustness of the convergence estimates which 
scale optimally with respect to the number of correctly executed subproblem solves.

In subsection \ref{sec14} we consider the practically more relevant situation of a distributed implementation on a network with predominantly local communication between unreliable compute cores. We slightly modify the approach taken in \cite{CXZ2017}, where a small amount of redundant storage capacity in the compute nodes is employed for a recovery of information lost due to hard faults, and apply the convergence theory of section \ref{sec1}. In particular,
we show theoretically and numerically that increasing the local redundancy improves the convergence behavior. We also briefly mention a server-client compute model from \cite{RiMo2017} composed of a reliable (fault-resilient) server network to which unreliable clients responsible for solving subproblems are attached, and demonstrate similar results. 

Let us uive some comments on the growing body of ABFT-related work. Even though it is
currently not clear which faults will be the most dominant and threatening in future HPC applications (in particular on exascale architectures), how to adequately model them, and which fault tolerance techniques will eventually produce the most impact in practice, the need in further research on fault tolerance techniques is undisputed. As to ABFT, core numerical linear algebra algorithms have been analyzed to a certain extent, and more recently various attempts have been made to address PDE solvers. We mention \cite{AGZ2015,AG2017a,AG2017b,AlG2017,ChDo2008,HGRW2016,RiMo2017,SaRi2017,St2017,StWe2015}, and refer to these papers for further references.

\subsection{Example: Overlapping domain decomposition}\label{sec12}
We present a simple variant of the overlapping DD method, and refer to \cite[Chapter 3]{TW2005} for more details. 
Consider a conforming FE discretization space $V\subset H_0^1(\Omega)$ of dimension $N=\dim(V)$ for the problem (\ref{MP}) on a quasi-uniform partition $T$ of $\Omega$ into cells of diameter $\approx h$.  Suppose that $\Omega=\cup_{i=1}^n \Omega_i$ is covered by a finite number of $n$ well-shaped subdomains $\Omega_i$ of diameter $\approx h_0$ which locally overlap. It is silently assumed that $h<<h_0$. Under natural assumptions on the alignment of the underlying FE partition $T$ with the boundaries 
of $\Omega$ and the $\Omega_i$, and for sufficient overlap of neighboring subdomains, a space splitting of the form
\be\label{DDS}
V = V_0 + \sum_{i=1}^n V_i,\qquad V_i:=V\cap H_0^1(\Omega_i), \quad i=1,\ldots,n,
\ee
is well-conditioned, with $\lambda_{\max}$, $\lambda_{\min}$, and $\kappa$ depending on the ellipticity constants $c_0,C_0$ of (\ref{MP}), the shape regularity of
the FE partition $T$ and the $\Omega_i$, and the overlap parameter $\delta$, but not on $h$ and $h_0$. To this end, (\ref{DDS}) must include   a properly constructed so-called coarse grid space $V_0$.  For $i=1,\ldots,n$ the operators $R_i: \, V_i\to V$ are the natural extension-by-zero operators
(the operator $R_0$ is special, an example is provided below). The bilinear forms $a_i(\cdot,\cdot)$ are inherited from $a(\cdot,\cdot)$ by restriction, i.e.,
$a_i(v_i,w_i)=a(R_iv_i,R_iw_i)$, $v_i,w_i\in V_i$, $i=0,1,\ldots,n$.
Neglecting the coarse problem associated with $V_0$ would result in a dependency of the splitting condition number $\kappa$
on roughly a factor $h_0^{-2}$. Note that there are more sophisticated space splittings
of DD type such as the Bank-Holst paradigm \cite{BH2003},
where the coarse problem is formally avoided by including a copy of it
into each of the subdomain problems with $i=1,\ldots,n$,  which we ignore here.

Within this framework, using the standard FE nodal basis in $V$ for representing elements of $V$ and $V_i\subset V$, $i=1,\ldots,n$, the matrix representations of the linear systems associated with subproblems on $V_i$ are given by smaller overlapping 
block-submatrices $A_i$ of the sparse matrix representation $A$ of the variational problem (\ref{VP})
associated with (\ref{MP}). This implies that the amount of data  
that needs to be communicated to enable a subproblem solve can be kept reasonable, i.e., linear in the dimension $M_i = \dim(V_i)$ of $V_i$, 
$i=1,\ldots,n$. If the partitioning is such that the number of locally overlapping $\Omega_i$ is small then communication is necessary only between a few processing units associated with subproblem solves for neighboring $\Omega_i$. In what follows the silent assumption is that subproblem solves for different $V_i$ take approximately the same time. Simplifying a bit, we therefore assume $M_i\approx M$, $i=1,\ldots,n$. Below, the constant $M$ is used to characterize the amount of storage and computational work per subproblem solve for $i=1,\ldots,n$. 
The coarse problem is typically generated by Galerkin discretization of (\ref{MP}) using a low-order FE space $V_0$ on a coarse partition $T_0$
with cell size comparable to the size of the subdomains $\Omega_i$. Since the coarse problem is in general treated differently, we do not specify $M_0=\dim(V_0)$. It is, however, clear that $M_0$ will scale linearly in $n$.

To discuss this in more detail, we make a number of simplifying assumptions concerning $V$ and the $V_i$, $i=0,1,\ldots,n$. Let the FE partition $T$ underlying $V$ be the refinement of a quasi-uniform coarse partition $T_0$ of element size 
$\approx h_0$.  Furthermore, assume that the cells of $T_0$ can be organized into $n$ disjoint groups of a few cells such that the union of the cells in each group represents a subdomain ${\Omega}'_i$ in a non-overlapping partition of $\Omega$ which is simply-connected and nicely shaped, and has only a small non-empty common boundary $\Gamma_{ii'}$ with a small number of direct neighbors ${\Omega}'_{i'}$, $i\neq i'$. Moreover, the subdomains  ${\Omega}'_i$  should be such  that the $V_i$ associated with their extensions $\Omega_i\supset \Omega'_i$ will have comparable dimensions $M_i\approx M$, $i=1,\ldots,n$. 

\begin{figure}[ht]
\centering
\includegraphics[width=0.48\textwidth,height=0.48\textwidth]{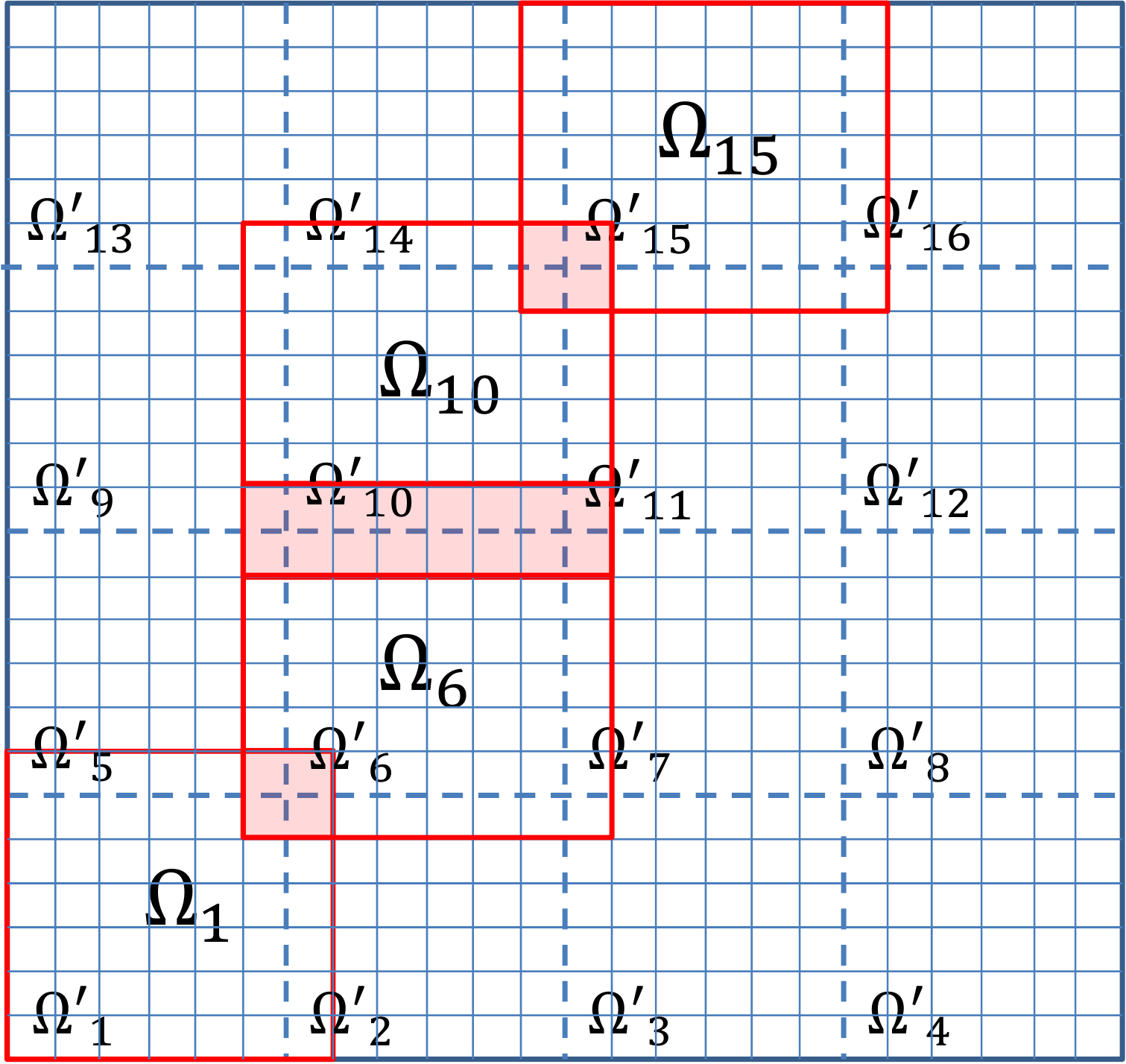}
\caption{Example of an overlapping domain partition for $d=2$ with $n=16$ subdomains $\Omega_i$, obtained from the squares $\Omega'_i$ of the coarse partition $T_0$ by adding one cell layer from the fine partition $T$. The parameters used for the figure
are $n_0=4$, $n_1=24$ ($k=6$), and $\ell=1$ ($\delta=1/6$). Only few of the domains $\Omega_i$ are depicted, some overlap regions $\Omega_i\cap \Omega_{i'}$ are highlighted.}
\label{Fig00}
\end{figure}

For the model situation of uniformly refined cube partitions of $\Omega=[0,1]^d$, 
choose some integers $n_0>1$, $k>1$, and $n_1=kn_0$, and
let $h_0=1/n_0$ and $h=1/n_1(=h_0/k)$ denote the mesh-widths of the cube partition $T_0$ and $T$,
respectively.  Then we can choose the $n=h_0^{-d}=n_0^d$ cubes in $T_0$ as subdomains ${\Omega}'_i$ which satisfy 
the above requirements: Each ${\Omega}'_i$, $i=1,\ldots,n$, has at most $3^d-1$ direct neighbors (for $d=3$ with face, edge, or vertex interfaces $\Gamma_{ii'}$ containing $\approx k^2$, $\approx k$, $1$ nodal points of $T$, respectively), and contains $\approx k^d$ nodal points from $T$, i.e., we have $M_i\approx M=k^d$.  Figure \ref{Fig00}
illustrates this construction for $d=2$.

This given, to achieve a mesh-independent condition number for the overall space splitting (\ref{DDS}), we can choose as $V_0\subset V$ a suitable FE space on $T_0$, e.g., a linear FE space on $T_0$ will do for a second-order elliptic problem such as (\ref{MP}). The FE spaces $V_i=V\cap H_0^1(\Omega_i)$, $i=1,\ldots,n$, depend on the choice of the larger subdomains $\Omega_i\supset \Omega'_i$. The traditional
overlapping scheme would form $\Omega_i$ as the union of all cells from $T$ in distance $\le \delta h_0$ from ${\Omega}'_i$,
where the overlap parameter  $\delta$ is a fixed positive number from $(0,1)$. Mild shape regularity assumptions on the resulting $\Omega_i$ will then guarantee robust condition number estimates of the form  
$$
\kappa\le C(1+\delta^{-1}),
$$
see \cite[Theorem 3.13]{TW2005}. Dropping the coarse grid problem, i.e., considering a space splitting of $V$ as in (\ref{DDS}) but without $V_0$,
would lead to the worse bound $\kappa\le Ch_0^{-2}(1+\delta^{-1})$.
Note that even though these estimates imply a deterioration of condition numbers proportional to $\delta^{-1}$ if $\delta\to 0$, in practice good performance has already been observed when $\Omega_i$ was obtained from $\Omega'_i$ by adding only a few layers of cells from $T$ around $\Omega'_i$. This may result in significantly smaller overlap regions $\Omega_i\cap\Omega_{i'}$ and smaller dimensions $M_i$ for the resulting $V_i$. 


\subsection{Master-slave network}\label{sec13}
We start with an idealized setting of a compute system consisting of a  reliable server $\mathcal{S}_0$ called master node with enough storage capacity to safely keep precomputed static arrays and master copies of dynamic data arrays, and a fixed number $p$ of unreliable compute nodes $\mathcal{C}_j$, $j=1,\ldots,p$, called slave nodes. 
We assume that slave nodes communicate with the master node but not with each other. During an iteration step, each slave node $\mathcal{C}_j$ is supposed to receive data and execution instructions to deal with a randomly assigned $V_i$ subproblem solve, $i=0,1,\ldots,n$,  and to return subproblem solutions to the master node.
All other work, such as forming linear combinations of vectors needed in 
(\ref{Rec0}) and (\ref{u}-\ref{v}), respectively, is performed by the master node $\mathcal{S}_0$.
Here, we allow for any value $p\le n+1$ which decouples the number of processors in the compute system from the size of the domain splitting. Note that we treat the $V_0$ subproblem in the same way as all other $V_i$ subproblems, $i=1,\ldots,n$, i.e., also 
the $V_0$ subproblem gets assigned to one of the slave nodes.
\begin{figure}[ht]
\centering
\includegraphics[width=0.60\textwidth,height=0.50\textwidth]{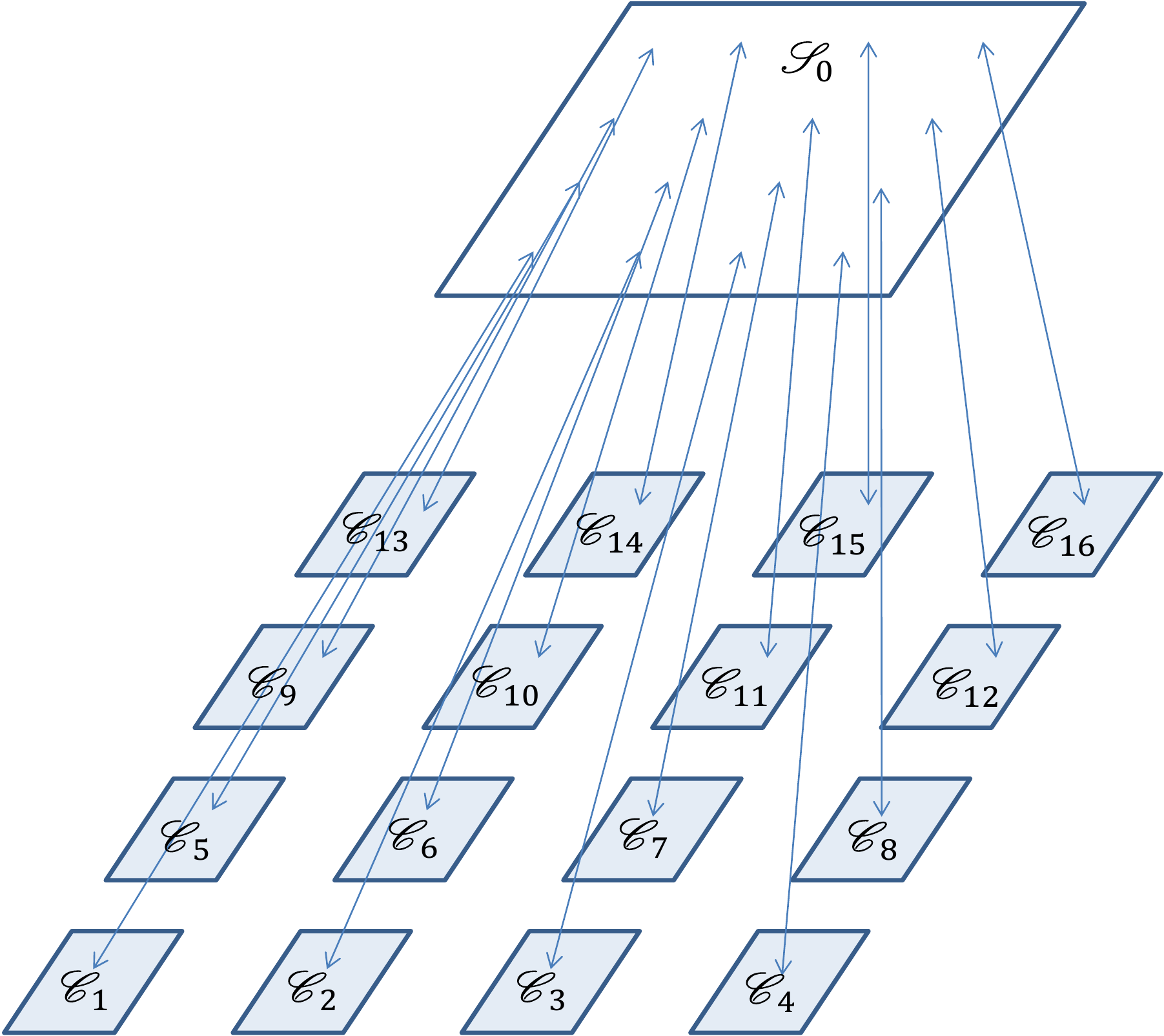}
\caption{Schematic view of the master-slave network. The compute nodes $\mathcal{C}_i$
communicate data directly to the reliable server $\mathcal{S}_0$ but not to each other. The server $\mathcal{S}_0$ needs to possess
storage capacity and compute power for global data arrays.}
\label{Fig0}
\end{figure}
 
Concerning the reason for and the nature of faults, we assume that faults are  detectable and represent unreturned or as wrong declared subproblem solves, i.e., we ignore soft errors such as bit flips in floating point numbers even if they were detectable. Moreover, the occurrence of a fault is unrelated to load balancing issues, i.e., slightly longer execution or communication times for a particular subproblem solve do not increase the chance of declaring such a process as faulty.  If the accelerated version (\ref{u}-\ref{v}) of the stochastic subspace correction scheme governed by Theorem \ref{theo2} is applied, then the number $f_m$ of faulty subproblem solves in iteration step $m$ is assumed a random variable 
with expectation $f$ and relatively small variance, independently of $m$ (this latter assumption is not needed for the iteration (\ref{Rec0})). It does not matter if faults are due to slave node crashes or communication failures, nor do we pose any restrictions on spatial patterns (which and how many slaves fail) or temporal correlations of faults (distribution of starting points and idle times
of failing slave nodes). 

Meaningful convergence results under such a weak fault model follow almost directly from our results in section \ref{sec1}, as long as we can select, uniformly at random and independently of previous iteration steps, $p$ subproblems out of the $n+1$ available ones  at the start of each iteration
step, assign them in a one-to-one fashion to the $p$ slave nodes,
and send the necessary data and instructions for processing the assigned subproblem solve to each of the $p$ slave nodes.
Indeed, if the time available for a solve step is tuned such that there is no correlation between faults and individual subproblem solves,
then one can safely assume that the index set $I_m$ corresponding to the $p_m=p-f_m$ as non-faulty detected subproblem solutions received by the master node by the end of the cycle is still a uniformly at random chosen subset of $\{0,1,\ldots,n\}$ that is independent of the index sets $I_0,\ldots,I_{m-1}$ used in the updates of the previous iteration steps.
It is important to realize that the latter independence property is the consequence of our scheme of randomly assigning subproblems to slave nodes, and not an assumption on the fault model.

Thus, Theorem \ref{theo1} applies, and yields according to
(\ref{EC}) the estimate 
\be\label{ECm}
\mathbb{E}(\|e_u^{(m)}\|^2) \le \prod_{s=0}^{m-1}\left(1-\frac{p_s}{\kappa (n+1)}\right) \|u\|^2,\qquad m=1,2,\ldots,
\ee
for the expected squared error if we formally set $\xi={\lambda}^{-1}_{\max}$. In practice, the value of $\xi$ can be determined by the steepest descent rule or from upper bounds $\bar{\lambda}$ for $\lambda_{\max}$. 

Similarly, Theorem \ref{theo2} gives guarantees for the expected squared error decay if we have safe a priori upper bounds $\bar{f}_m\ge f_m$ for the number of faults and bounds $\underline{\lambda}\le \lambda_{\min}\le \lambda_{\max}\le\bar{\lambda}$ for the spectrum of $P$.   
In this case the subproblem solves  can be performed using parameters $\alpha_m,\beta_m$ determined from (\ref{AM}), (\ref{BM}) with $p_m=p-f_m$ replaced by $\underline{p}_m=p-\bar{f}_m$, and (\ref{ECA}) implies
\be\label{ECAm}
\mathbb{E}(\|e_u^{(m)}\|^2) \le 2\prod_{s=0}^{m-1}\left(1-\frac{\underline{p}_s}{\sqrt{\bar{\kappa}} (n+1)}\right) \|u\|^2,\qquad m=1,2,\ldots,
\ee
where $\bar{\kappa}=\bar{\lambda}/\underline{\lambda}\ge \kappa$, compare the derivation of (\ref{ECAU}). 

\medskip
We have conducted some preliminary numerical experiments using the example of the model problem (\ref{MP}) with homogeneous diffusion coefficient $a(x)=1$ and right-hand side $f(x)=1$ for $d=2$.  
The domain $\Omega$ is the unit square, equipped with a uniform coarse square partition $T_0$ of step-size $h_0=1/n_0$ in each direction,
and a uniform fine square partition $T$ of step-size $h=1/{n_1}$, where $n_0$ divides $n_1$
(in other words, $n_1=kn_0$ for some integer $k$). 
Both $V$ and $V_0$ are given by
bilinear finite element spaces on the respective square partitions. The stiffness matrix  $A$ and right-hand side $b$
are computed exactly. 
The subdomains $\Omega_i$, $i=1,\ldots,n$, are obtained by adding to each of the squares in $T_0$ in each direction  $\ell$ layers of  
square cells from the fine partition, i.e., the overlap parameter is $\delta=\ell h/h_0=\ell/ k$. Below, we report numerical results for the values $n_0=20$, $n_1=400$, $\ell=6$ which gives $k=20$, $\delta=0.3$, and results in a overlapping partition of $\Omega$ with $n=400$ subdomains $\Omega_i$. The associated DD space splitting (\ref{DDS}) with weights $\omega_i=1$, $i=0,1,\ldots,n$, is well-conditioned, with a value $\kappa\approx 6$. All subproblems, including the $V_0$ subproblem, have approximately the same dimension $M_i\approx M=400$, $i=0,1,\ldots,n$. Despite the fact that the dimension $N\approx 160000$ of the discretization space $V$ is still moderate, our numerical findings for this parameter set can be considered as sufficiently representative. 
 
Iterations are always started from the zero vector, and terminated when a relative error reduction of $\epsilon_0=10^{-6}$ is achieved. Here,
errors are given by the error indicators $\epsilon$ computed as explained at the end of the appendix.
For simplicity, all subproblems are solved by a sparse elimination method.  
In the numerical experiments with the one-step iteration (\ref{Rec0}), for the relaxation parameters $\xi_{m,i}=\xi_m \omega_i$ the value $\xi_m$ was determined  by the steepest descent rule (\ref{XiOpt}). As a matter of fact, and in agreement with the remarks after Theorem 1, iteration counts with constant values $\xi_m=\xi$ were typically
higher, and are not reported here. For the accelerated iteration, the constants $\xi_m=\bar{\lambda}^{-1}$ and $\eta_m=(\underline{\lambda}\bar{\lambda})^{-1/2}$ in (\ref{u}-\ref{v}) have been obtained after initial tests with values $\bar{\lambda}=3.33$ and
$\underline{\lambda}=0.9$, such that iteration counts were near-optimal for the standard deterministic iteration, where $p_m=n+1$ and $I_m=\{0,1,\ldots,n\}$. Note that $\lambda_{\max}\approx 5$ for our problem and thus this choice of $\bar{\lambda}$ corresponds to overrelaxation.

\begin{figure}[ht] 
\centering
\includegraphics[width=0.46\textwidth,height=0.26\textheight]{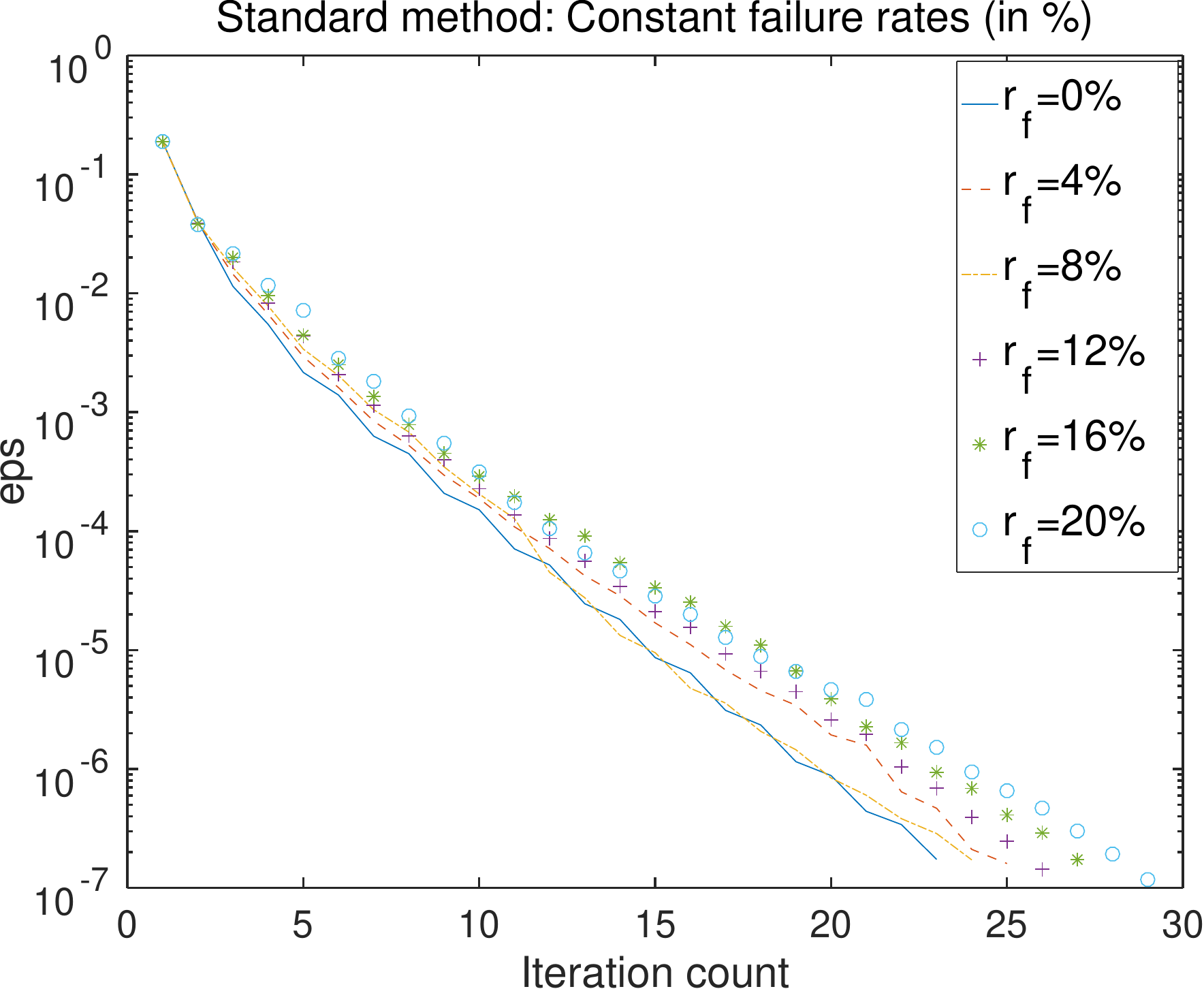} $~~~~~~~~~~~$
\includegraphics[width=0.46\textwidth,height=0.26\textheight]{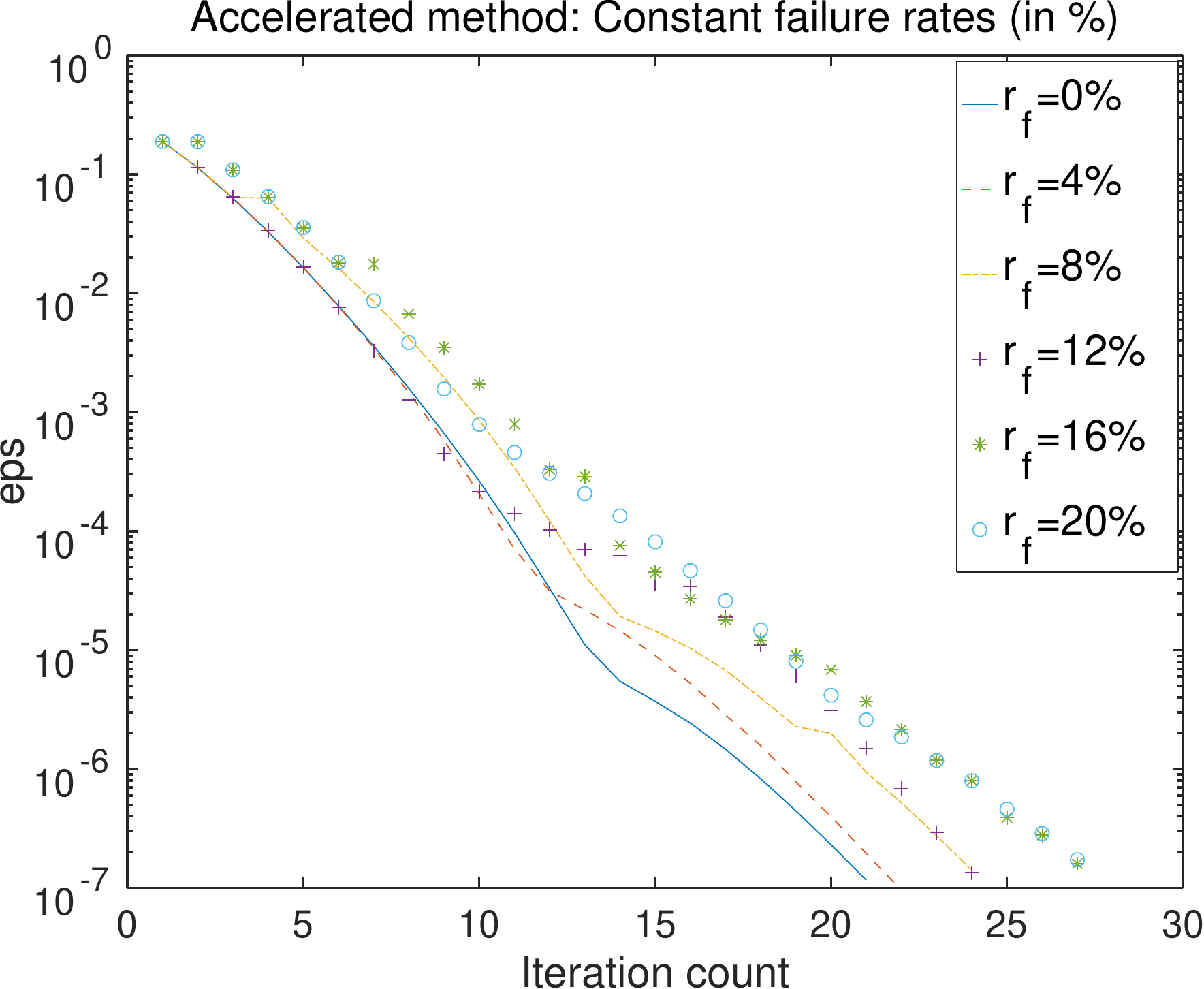}
\caption{Convergence results for a model DD space splitting on the master-slave network with different constant failure rates. Left:
Results for the one-step method (\ref{Rec0}) with steepest descent choice (\ref{XiOpt}) for $\xi_m$. Right: Results for the accelerated method (\ref{u}-\ref{v}) with near-optimal parameters $\xi$ and $\eta$.}
\label{Fig2}
\end{figure}

In Figure \ref{Fig2} we show convergence results for the above described model problem for different constant failure rates $r_f$, where for simplicity
the number of available compute nodes is set to $p=n+1$. A constant failure rate $r_f\in [0,1]$ means that throughout the recursions we have
chosen for $p_m$ the constant value $p^\ast=\lfloor (1-r_f)(n+1)\rfloor$. Thus, in each step the fixed number $f^\ast=n+1-p^\ast$ of compute nodes 
fails to return correct subproblem solutions. Then, in agreement with our assumption of a random assignment of subproblems to compute nodes in
each iteration step, the index set $I_m$ was selected as a random subset of size $p^\ast$ from $\{0,1,\ldots,n\}$. 
For each value $r_f$ we show plots of error indicators $\epsilon$ as function of $m$ for one run only, i.e., one sequence of index sets $I_m$, and not the expectation $\mathbb{E}(\|e^{(m)}\|)$ of the true error. Due to further experiments, this seems fully justified since error indicator curves
for different runs were only marginally different,
and almost identical with those for the true error $\|e^{(m)}\|$. Moreover, testing of more realistic fault scenarios, where the $f_m$ were generated by an independent
sequence of random integers uniformly distributed in a certain interval $[f^\ast - \Delta f,f^\ast +\Delta f]$ yielding the same failure rate
$r_f$ but only in expectation, revealed that convergence behavior and iteration counts are very robust to the variance of randomly created sequences
$f_m$ and $p_m=n+1-f_m$, respectively,  as long as $r_f:=\mathbb{E}(f_m/(n+1))$ was fixed. In other words, we consider Figure \ref{Fig2} as a fair illustration of the convergence properties of our fault mitigation approach for an overlapping DD method on a master-slave compute network as proposed in this subsection. 

Figure \ref{Fig2} (left) shows the decay of the error indicator $\epsilon$ for the one-step method (\ref{Rec0}) with $\xi_m$ determined by the steepest descent rule (\ref{XiOpt}) for different constant failure rates $r_f\in [0,0.2]$.
The case $r_f=0$ (bold line) corresponds to a non-faulty compute network, where all $p=n+1$ processors return valid subproblem solutions ($I_m=\{0,1,\ldots,n\}$). For our particular test problem, the required relative error reduction 
of $\epsilon_0=10^{-6}$ was reached after $23$ iteration steps.
When the failure rate $r_f$ is increased, the number of iteration steps to termination slightly grows, see Table \ref{tab0}, which is visible from the error curves as well. This is in agreement with Theorem \ref{theo1} (and the comments following it) which predicts an upper bound for the expected error reduction per iteration step of at least
$(1-(1-r_f)/\kappa)^{1/2}$ since we have
$$
\mathbb{E}(\|e^{(m+1)}\|^2) \le (1-\frac{p^\ast}{(n+1)\kappa})
\mathbb{E}(\|e^{(m)}\|^2)=(1-\frac{1-r_f}{\kappa})\mathbb{E}(\|e^{(m)}\|^2),\qquad m=0,1,\ldots.
$$
Recall that the error reduction per step for (\ref{Rec0}) with the steepest descent rule (\ref{XiOpt}) is at least as good as with any fixed choice $\xi_m=\xi$ for the relaxation parameter.

\begin{table}[ht]
\centering
\begin{tabular}{l|rrrrrr}
&\multicolumn{6}{c}{Iteration counts for different $r_f$}\\
\hline
\;\; Method & $r_f=0$ & $r_f=0.04$ & $r_f=0.08$ & $r_f=0.12$& $r_f=0.16$& $r_f=0.2$\\
\hline
(\ref{Rec0}), steepest descent (\ref{XiOpt}) & $23\;$ & $25\;$ &$24\;$ &$26\;$ & $27\;$&$29\;$ \\
(\ref{Rec0}), $\xi_m=0.4$ &  $29\;$ &$30\;$ &$30\;$ & $31\;$&
$33\;$ &$40\;$\\
(\ref{u}-\ref{v}), $\xi_m=0.3$, $\eta_m=0.58$ &  $21\;$ &$22\;$ &$24\;$ & $24\;$&
$27\;$ &$27\;$
\end{tabular}
\caption{Iteration counts for reaching a relative error reduction of $\epsilon_0=10^{-6}$ for the model test problem and the iteration (\ref{Rec0}) with $\xi_m$ determined by the steepest descent rule (\ref{XiOpt}) and with near-optimal constant $\xi_m=0.4$, and the accelerated iteration
(\ref{u}-\ref{v}) with near-optimal constant values
$\xi_m=0.3$ and $\eta_m=0.577$.}
\label{tab0}
\end{table}
In Figure \ref{Fig2} (right), we show similar results for the 
accelerated method (\ref{u}-\ref{v}), see also Table \ref{tab0} for the recorded iteration counts to termination. The parameters $\xi=0.3$ and $\eta=0.577$ were determined by experiment, and are near-optimal in the sense that for them the iteration count of the additive Schwarz method for the given problem and error reduction level $\epsilon_0$ is close to minimal. Again, the graphs show that
our method behaves according to theory. That the iteration counts are almost the same as for (\ref{Rec0}) is at first glance unexpected but can be explained as follows. On the one hand, for well-conditioned space splittings with $\kappa<10$,
such as the overlapping  DD space splittings for our test problem, one should not expect dramatic gains. On the other hand,
the accelerated method is run with constant, although near-optimal relaxation parameters while the steepest descent rule for (\ref{Rec0}) results in a nonlinear iterative scheme which is superior to any iteration (\ref{Rec0}) with constant $\xi_m=\xi$. If one compares the accelerated method with any of the latter methods, the improvement by acceleration becomes more visible. This is supported by the iteration counts to termination for (\ref{Rec0}) with constant $\xi_m=\xi=0.4$ reported in 
Table \ref{tab0} (this value is near-optimal in the sense 
described above). Still, an application of the accelerated scheme in the case of well-conditioned space splittings such as the overlapping domain DD scheme for (\ref{MP}) is questionable, as the possible gains do not justify the additional effort needed for parameter tuning and the iteration itself.

We want to stress that an implementation of the above  random assignment scheme
on a master-slave compute network is only of academic
interest as it neglects the massive communication overhead necessary before and after each compute cycle. Recall that the advantage of the random assignment scheme is that it enforces almost automatically, under very mild requirements for the fault model, the independence assumption for the random index sets $I_m$ of subproblem solves actually used in each iteration step that is needed in the proofs of 
Theorem \ref{theo1} and \ref{theo2}. The obtained convergence estimates (\ref{ECm}) and (\ref{ECAm}) are in some sense the best possible ones since they signal a loss of convergence speed compared to a fault-less environment which is only proportional to the fault rate $f_m/p$.

Parallelization gains for the accelerated method are also questionable since update steps have to be executed on the master node  to keep dynamic data correct (for a more detailed estimate of the
runtime of a parallel implementation of the DD method on this and other compute networks discussed below we refer to subsection \ref{sec15}).
In particular, the accelerated iteration (\ref{u}-\ref{v}) needs additional global vector operations which seem prohibitive, especially if
$p<<n+1$, i.e., if in each iteration step only a small percentage of the subproblems can be assigned to
the slave nodes. Remedies for this problem have been discussed in the literature for similar problems,
see \cite{FerRic2016,LeeSid2013}.  

In the next subsection, we will try to achieve a better compromise between communication overhead and overall parallel efficiency  on the one hand, and a matching of the  theoretical assumptions for our theory outlined in section \ref{sec1} on the other hand. Needless to say that we will not be able to achieve a fully optimal 
solution (in particular with regard to solving the coarse grid problem and computing global quantities which represents a challenging bottleneck
for parallelization), and that further assumptions on the fault model may become necessary. The major 
difference will be that we give up the random assignment scheme of subproblem solves to compute nodes, and accept the common paradigm 
of distributed data storage and program execution paired with maximally local communication that is typical
for parallelization efforts for DD and multigrid methods developed during the last 30 years.

\subsection{Local communication network}\label{sec14} 
Standard parallelization efforts for DD and multigrid algorithms
under a no-fault assumption avoid global storage and global communication steps as much as possible, compare, e.g., \cite{BH2003}. In a typical parallel implementation, the assignment of subproblems to compute nodes is fixed,
and data storage is decentralized. When it now comes to compute node failures, one is additionally confronted with the possibility of the loss of dynamic and static data arrays
associated with one or several subproblems. This is in contrast to what we have assumed for the previously discussed master-slave network, where safe copies of all data arrays were maintained at the master at the expense of high communication cost. 
Consequently, one now needs to incorporate local redundancy in data storage and/or restoration methods for missing data. Moreover, the localization in time and space (i.e., with regard to the neighbor structure of the domain partition into subdomains of faulty compute nodes becomes an issue because it may contradict the randomness and independence assumption {\bf A} for the index sets $I_m$ corresponding to the correctly executed subproblem solves in each iteration step.

In the literature, several proposals already exist using different settings. 
The majority of papers ignores the loss of static data arrays
associated with a subproblem, and concentrates on restoring dynamic data arrays that change during the iteration. Most of the existing analysis is on faults isolated in time and space. Often, the time to restart 
a failing compute node is ignored.  For example, \cite{HGRW2016} uses the concept of a so-called local superman unit (i.e., additional, more powerful compute nodes) that can be setup in no time to replace a failed compute node, and has enough compute power and speed to recover lost local solutions from scratch much faster than normal compute nodes. Alternatively, \cite{CXZ2017} proposes the introduction of redundancy
in local storage such that when a compute node executing a given $V_i$ subproblem fails, there is a neighboring node that has all current data arrays for the $i$-th subproblem. Thus, recovery does not need to start from scratch. Several strategies concerning which compute node performs which recovery action during the time the failed compute node is 
not back are discussed in \cite{CXZ2017} yielding different deterministic error decay guarantees. Both papers assume that faults are spatially 
isolated, i.e., neighbors of a failing compute node that keep copies of its data arrays do not fail at the same time.

\begin{figure}[ht]
\centering
\includegraphics[width=0.65\textwidth,height=0.57\textwidth]{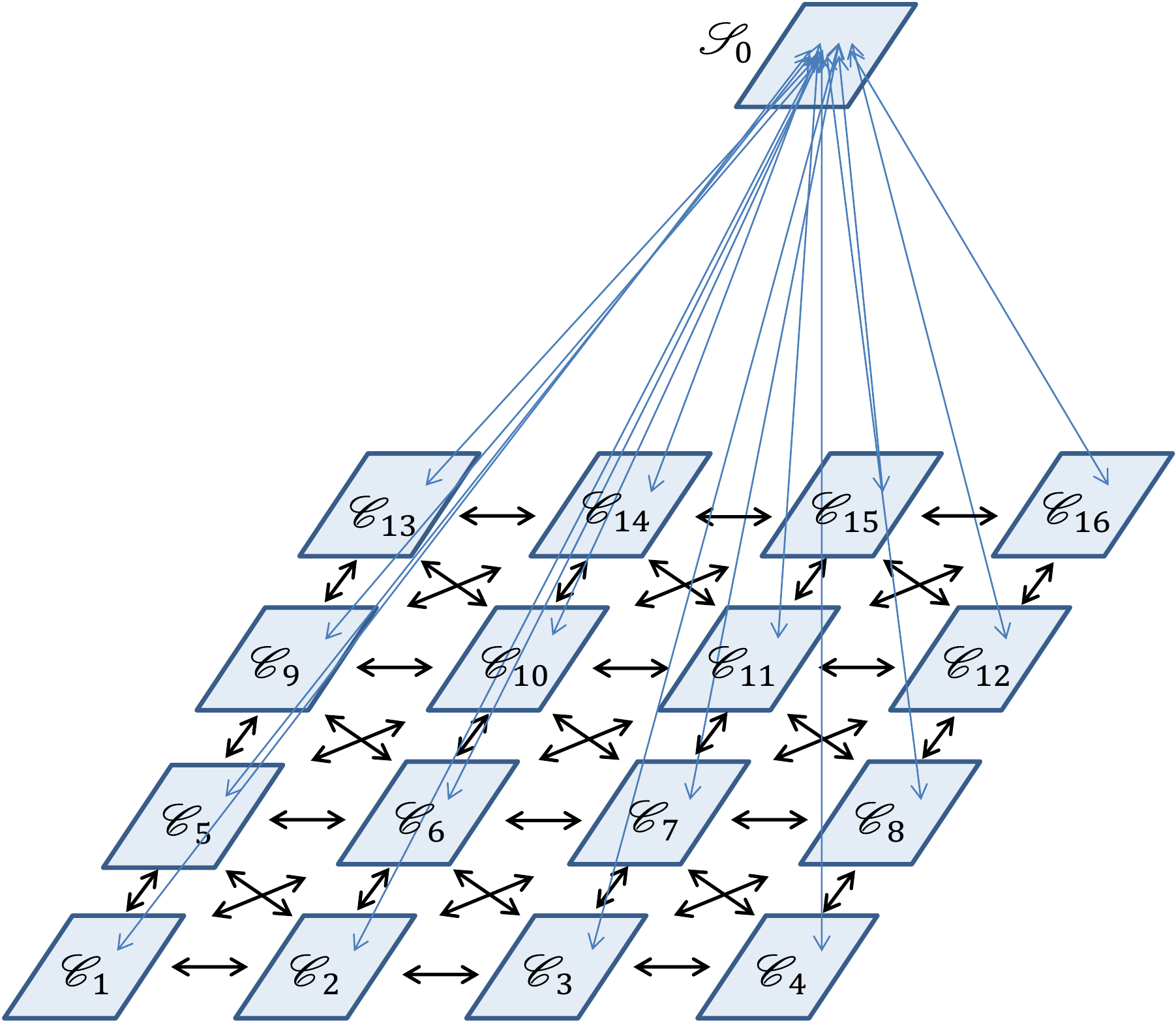}
\caption{Schematic view of the local communication network. High load communication
between the compute nodes $\mathcal{C}_i$ is local while the communication between
compute nodes $\mathcal{C}_i$ and the reliable server $\mathcal{S}_0$ concerns only small data arrays. In contrast to the master-slave network, compute power and storage capacity of $\mathcal{S}_0$ and the compute nodes $\mathcal{C}_i$ can be of the same order. }
\label{Fig1}
\end{figure}
We incorporate these ideas in slightly modified form into our discussion of possible applications of the convergence results of section \ref{sec1}. Consider a network of (at least)
$n$ unreliable compute nodes $\mathcal{C}_i$, $i=1,\ldots,n$, and a reliable server $\mathcal{S}_0$, i.e., we have $p=n+1$ for the number of available processors.  The subproblems associated with the DD space splitting (\ref{DDS}) will be statically assigned such that data arrays for the $V_i$ subproblem are stored at $\mathcal{C}_i$, $i=1,\ldots,n$. Unless a failure occurs, the main task of $\mathcal{C}_i$ is then to solve
the $V_i$ subproblem in each iteration step. Moreover, the bottleneck problem associated with $V_0$ and other global quantities such as error  estimators and scalar products is assigned to the reliable server $\mathcal{S}_0$ (to achieve reliability, one can apply standard fault tolerance techniques such as node replication). In contrast to the master-slave network
considered in subsection \ref{sec13}, where $\mathcal{S}_0$ had basically a storage function and performed the update steps but not the coarse problem solve,
now $\mathcal{S}_0$ only needs to receive from and communicate back to all other compute nodes small data arrays before and after each compute cycle, respectively, to solve the $V_0$ subproblem and to compute other global quantities (see the appendix for more details). Communication of larger data arrays between compute nodes $\mathcal{C}_i$ (typically dynamic data arrays but possibly also  static data arrays associated with a $\Omega_i$ in case of failure of the associated compute node $\mathcal{C}_i$)  is local, i.e.,
only compute nodes $\mathcal{C}_{i}$ and $\mathcal{C}_{i'}$ responsible for  $V_i$ and $V_{i'}$ subproblem solves with overlapping $\Omega_i\cap \Omega_{i'}\neq \emptyset$  need to communicate. Moreover, since the number of neighbors is uniformly bounded by a certain constant
\be\label{Locall}
\bar{l}\ge \max_{i=1,\ldots,n}\,|\{i': \Omega_i\cap \Omega_{i'}\}|
\ee 
(for the example of the DD space splitting induced by cube partition discussed in subsection \ref{sec12}, we have $\bar{l}=3^d-1$), this 
communication can be implemented in parallel for all compute nodes $\mathcal{C}_i$ in a finite number of sweeps proportional to $\bar{l}$.

This  is, up to the details of dealing with the coarse problem, the standard approach to parallelizing DD and multigrid methods. We keep
the option of adding additional compute nodes on the fly to replace processors that fail over a period of many cycles or permanently, thus rejuvenating the whole system (this is analogous to the superman unit concept from \cite{HGRW2016}). If such an additional compute node takes over a particular $V_i$ solve, connections with all nodes responsible for $V_{i'}$ solves with $\Omega_i\cap \Omega_{i'}\neq \emptyset$ need to be enabled in an update of the communication structure. We also borrow from \cite{CXZ2017} the idea of keeping, at each node, redundant copies of the data arrays
of a few neighboring compute nodes (as long as the copies come from direct neighbors with $\Omega_i\cap \Omega_{i'}\neq \emptyset$, this does not require significant changes in the communication structure). In contrast to \cite{CXZ2017}, where always two compute nodes are grouped in pairs and both keep the data arrays of a pair of subproblem solves, we assume that each set of dynamic and static data arrays assigned to $\mathcal{C_i}$ has up-to-date copies in at least $l\ge 1$ neighboring compute nodes $\mathcal{C}_{i'}$ ($l = 1$ corresponds to the analysis in 
\cite{CXZ2017}). Although larger values $l$ lead to a proportional increase of storage and communication overhead at the compute nodes, they decrease at the same time the chance of complete loss of data associated with a subproblem and also give some flexibility of artificially enforcing the randomness and independence assumptions for the index sets $I_m$ that are prerequisite for our convergence proofs. 

We now describe our methodology of dealing with faults. We start with the assumption of spatially and temporarily isolated faults, i.e., if after a given cycle the subproblem solve for $V_i$ is detected as faulty and its associated compute node is not available for a certain number of cycles, then all neighboring compute nodes remain non-faulty during the whole time. That a compute note will be unavailable for a very long time or permanently can be counteracted by adding compute nodes to the network. In the update step following the fault detection, we proceed as usual,
however, no communication to the failing node is possible, and the redundant dynamic data arrays for $V_i$, that are stored at neighboring compute nodes, are updated assuming zero change coming from the failing node. In subsequent compute cycles, the $\ge l$ neighboring compute nodes, that keep copies of 
the $V_i$ data, change their role temporarily as follows: For simplicity, select exactly $l$ such neighboring nodes, and denote their indices by $j_1,j_2,\ldots,j_l (\neq i)$. Choose, with equal probability $1/(l+1)$, an index from $\{i,j_1,\ldots,j_l\}$. If the chosen index is $i$ then each of the $l$ selected compute nodes solves its originally assigned subproblem (as a consequence, the $V_i$ subproblem solve is not executed in this cycle). If the chosen index is $j_s$ for some $s=1,\ldots,l$, then the compute node $\mathcal{C}_{j_s}$
assumes the role of the failing node, i.e., it executes a $V_i$ subproblem solve instead of its statically assigned $V_{j_s}$ subproblem solve.
The other $l-1$ selected compute nodes compute their originally assigned subproblems. Consequently, the $V_{j_s}$ subproblem is not executed in this cycle. If $l=1$, to avoid repetition, this rule can be modified as follows: The single selected neighboring compute node $\mathcal{C}_{j_1}$ that keeps redundant $V_i$ information will then alternately execute the $V_i$ and its own $V_{j_1}$ subproblem solve (this was the proposal in \cite{CXZ2017}). 

The proposed random assignment rule in the neighborhood of a failing node needs only
a certain amount of additional local communication, see the appendix for more details. 
We proceed with it until the failing node is restarted or replaced by a new node.  At this moment, all information (static, dynamic, and redundant) that was previously owned by the $V_i$ node 
needs to be recollected from the neighbors. Compared to the existing local communication between neighboring compute nodes, the communication overhead is now only in the static data. 

The assumption of fault locality (in space and time) is often made in the literature for analysis purposes. Multiple faults at the same time are not an
issue, as long as they remain spatially disjoint (i.e., the neighborhoods of different failing nodes do not intersect), they are covered by the 
convergence theory outlined below. It is clear that one can deal with neighbor pairs or even larger local groups of failing compute nodes by designing similar repair rules. The occurrence of large spatially correlated parts of the compute network is a situation that is beyond the scope of ABFT methods with restrictions on the amount of global communication during the algorithm execution. 
We will not further discuss this issue. As to the temporal distribution of faults, some researchers \cite{PAS2014} assume a Weibull distribution for the failure arrival times whereas little is known about the length of
failing. The Weibull distribution function is given
by
$$
F_W(t;k,\lambda)= \left\{\ba{ll} 1-e^{-(\lambda t)^k},&t\ge 0,\\ 0,&t<0,\ea  \right.
$$
and is characterized by two positive parameters, the scale parameter $\lambda$ and the shape parameter $k$. It generalizes the exponential distribution ($k=1$). We will use the Weibull distribution in our numerical simulations reported below.

We claim that minor modifications in the argument leading to Theorem \ref{theo1} allow us to formulate convergence results if the above strategy is 
followed, and the faults are spatially isolated and occur uniformly at random and independently in the compute network. 
We argue again for (\ref{Rec0}). At the iteration step from $u^{(m)}$ to $u^{(m+1)}$
the random set $I_m$ of correctly executed subproblem solves is given by
$$
I_m = \{0,1,\ldots,n\}\backslash \{i_1,\ldots,i_{f_m}\},\qquad i_1,\ldots,i_{f_m}\in \{1,\ldots,n\},
$$
where $f_m$ is the number of failing nodes during the cycle (these consist of newly detected faults, and nodes that failed in previous cycles but
are still under repair). 

Unfortunately, the $I_m$ are not uniformly at random selected subsets of $\{0,1,\ldots,n\}$. 
To this end, it is instructive to first look for a convergence estimate in the case of a single fault. Then $f_m=1$ and $I_m=\{0,1,\ldots,n\}\backslash \{i^\ast_m\}$ for $m=m_0,\ldots,m_1$, where for $m=m_0$ the $V_{i^\ast}$ subproblem solve was detected as faulty, and after the step with $m=m_1$ the failing compute node $\mathcal{C}_{i^\ast_m}$  was restarted. For $m = m_0$ the index $i^\ast_m$ equals  an integer $i^\ast$ uniformly at random selected from $\{1,\ldots,n\}$ (the index of the subproblem solve assigned to the faulty node), while for  $m = m_0+1,\ldots,m_1$ the index $i^\ast_m$ is selected uniformly at random from the index set 
\be\label{Jset}
J:=\{i^\ast,j_1,\ldots,j_l\}
\ee 
associated with the failing node. For all $m\not\in \{m_0,\ldots,m_1\}$, we have $I_m=\{0,\ldots,n\}$. It is not hard to check that these $I_m$ are a sequence of independent index sets whenever $1\le m_0\le m_1$ are fixed (this follows from the fault model and our above random assignment rule of selecting an index from $J$).

The iteration steps (\ref{Rec0}) with $m < m_0$
and $m > m_1$ are deterministic, and thus
\be\label{A11}
\|e_u^{(m+1)}\|^2\le (1-\frac{1}{\kappa})\|e_u^{(m)}\|^2 
\ee
if we set $p_m=n+1$ and $\xi=\lambda_{\max}^{-1}$ in Theorem \ref{theo1} for these $m$. For $m=m_0$, the failing node $i^\ast$ can be considered chosen uniformly at random from $\{1,\ldots,n\}$
(recall that the $V_0$ subproblem solve is treated by a the reliable server $S_0$ and is never faulty). Thus, $I_{m_0}=\{0\}\cup I'_{m_0}$ 
where $I'_{m_0}$ is a uniformly at random selected subset of  $\{1,\ldots,n\}$ of size $n-1$. Similarly, for $m= m_0+1,\ldots,m_1$ we have
$I_{m}=(\{0,1,\ldots,n\}\backslash J)\cup I'_{m_0}$, where now $I'_{m_0}\subset J$ is an index set of size $l$, of uniformly at random selected indices from $J=\{i^\ast,j_1,\ldots,j_l\}$ defined in (\ref{Jset}). This situation is not directly covered by Theorem \ref{theo1}. However, a small modification of its proof gives the following estimate
for the conditional expectation of the squared error in one recursion step (\ref{Rec0}) under slightly different assumptions.

\begin{cor}\label{cor1} Let the disjoint sets $I^s$, $s=1,\ldots,S$ form a partition of $\{0,1,\ldots,n\}$, i.e.,
$$
I^1\cup \ldots\cup I^S=\{0,1,\ldots,n\},\qquad I^s\cap I^{s'}=\emptyset,\quad s\neq s'.
$$
For fixed $u^{(m)}$, let  $u^{(m+1)}$ be given by (\ref{Rec0}) with $I_m=I_m^1\cup \ldots\cup I_m^S$, where the
$I^s_m$  are uniformly at random selected subsets of $I^s$ of
size $p^s\le |I^s|$ , $s=1,\ldots,S$.
Then setting 
$$
\bar{r}_p = \max_{s=1,\ldots,S} \frac{p^s}{|I^s|},\qquad
\underline{r}_p = \min_{s=1,\ldots,S} \frac{p^s}{|I^s|}
$$
and taking $\xi_m=\underline{r}_p/(\bar{r}_p\lambda_{\max})$ yields the estimate
\be\label{CorE}
\mathbb{E}(\|e_u^{(m+1)}\|^2|u^{(m)}) \le \left(1-\frac{\underline{r}_p^2}{\bar{r}_p \kappa}\right) \|e_u^{(m)}\|^2.
\ee
\end{cor}

{\bf Proof}. Indeed, compared to the proof of Theorem \ref{theo1} the only changes are in the evaluation of 
\bea
\mathbb{E}(a(e_u^{(m)},\sum_{i\in I_m} \omega_i R_iT_ie_u^{(m)})|u^{(m)})&=&\sum_{s=1}^S
a(e_u^{m}, \mathbb{E}(\sum_{i\in I^s_m} \omega_i R_iT_ie_u^{(m)}))\\
&=& \sum_{s=1}^S \frac{p^s}{|I^s|}a(e_u^{(m)},\sum_{i\in I^s}\omega_i R_iT_ie_u^{(m)})\\
& \ge & \underline{r}_p a(Pe_u^{(m)},e_u^{(m)}),
\eea
and of
\bea
\mathbb{E}(\|\sum_{i\in I_m} \omega_i R_iT_ie_u^{(m)}\|^2|u^{(m)})&\le&
\lambda_{\max} \mathbb{E}(\sum_{i\in I_m} \omega_i a_i(T_ie_u^{(m)},T_ie_u^{(m)})|u^{(m)})\\
&=&\lambda_{\max} \mathbb{E}(a(e_u^{(m)},\sum_{i\in I_m} \omega_i R_iT_ie_u^{(m)})|u^{(m)})\\
&\le& \lambda_{\max}\bar{r}_p a(Pe_u^{(m)},e_u^{(m)}).
\eea
Substituting these upper and lower estimates into the corresponding expressions of the formula for $\mathbb{E}(\|e_u^{(m+1)}\|^2|u^{(m)})$, we get
\bea
\mathbb{E}(\|e_u^{(m+1)}\|^2|u^{(m)}) &=& \mathbb{E}(\|e_u^{(m)}\|^2-2\xi
a(e_u^{(m)},\sum_{i\in I_m} \omega_i R_iT_ie_u^{(m)}) +\xi^2
\|\sum_{i\in I_m} \omega_i R_iT_ie_u^{(m)}\|^2|u^{(m)})\\
&=& \|e_u^{(m)}\|^2-2\xi
\mathbb{E}(a(e_u^{(m)},\sum_{i\in I_m} \omega_i R_iT_ie_u^{(m)})|u^{(m)}) \\
&&                    \qquad\qquad\qquad +\xi^2
\mathbb{E}(\|\sum_{i\in I_m} \omega_i R_iT_ie_u^{(m)}\|^2|u^{(m)})\\
&\le&
\|e_u^{(m)}\|^2 - 2\xi \underline{r}_p a(Pe_u^{(m)},e_u^{(m)})+\lambda_{\max}\bar{r}_p\xi^2 a(Pe_u^{(m)},e_u^{(m)})\\
&\le& (1-\frac{\lambda_{\max}\xi(2\underline{r}_p-\lambda_{\max}\bar{r}_p\xi)}{\kappa})\|e_u^{(m)}\|^2.
\eea
This bound is optimized for $\lambda_{\max}\xi=\underline{r}_p/\bar{r}_p$, and yields (\ref{CorE}). \hfill $\Box$

\medskip
For $S=1$ we recover the statement of Theorem \ref{theo1} since then $\underline{r}_p=\bar{r}_p=p_m/(n+1)$. 
Note that $\frac{p^s}{|I^s|}$ can be interpreted as the rate of non-faulty compute nodes $\mathcal{C}_i$ with indices in 
$I^s$, and that $\underline{r}_p$ and $\bar{r}_p$ stand for the minimum and maximum of these rates.

For $m=m_0$, we apply Corollary \ref{cor1} with $S=2$, $I^1=\{0\}$, $p^1=1$, $I^2=\{1,\ldots,n\}$, $p^2=n-1$. Then $\underline{r}_p=(n-1)/n$, $\bar{r}_p=1$, and we obtain
\be\label{A22}
\mathbb{E}\|e_u^{(m_0+1)}\|^2|u^{(m_0)}) \le \left(1-\frac{ (n-1)^2}{n^2 \kappa}\right) \|e_u^{(m_0)}\|^2,
\ee
if we set $\xi_{m_0}=(n-1)/(n\lambda_{\max})$. 
For the steps with $m = m_0+1,\ldots,m_1$ we set $S=2$, $I^1=J=\{i^\ast,j_1,\ldots,j_l\}$, $p^1=l$, $I^2=\{0,1,\ldots,n\}\backslash J$, $p^2=|I^2|$, and with $\underline{r}_p=l/(l+1)$, $\bar{r}_p=1$, and $\xi_m=l/((l+1)\lambda_{\max})$ in Corollary \ref{cor1}, this yields
\be\label{A33}
\mathbb{E}(\|e_u^{(m+1)}\|^2|u^{(m)}) \le \left(1-\frac{l^2}{(l+1)^2 \kappa}\right) \|e_u^{(m)}\|^2, \qquad m = m_0+1,\ldots,m_1.
\ee
for the conditional expectation of $\|e_u^{(m+1)}\|^2$.

Using the independence of the index sets $I_m$, according to (\ref{A11}), (\ref{A22}), (\ref{A33}) we get guaranteed one-step reduction factors for the (unconditional) expectation of the squared error $\mathbb{E}(\|e_u^{(m)}\|^2)$
of
$$
(1-1/\kappa),\qquad (1- ((n-1)/n)^2/\kappa),\qquad (1-(l/(l+1))^2/\kappa),
$$
for the different types of iteration steps corresponding to $m\not\in\{m_0,\ldots,m_1\}$, $m=m_0$, and $m=m_0 +1,\ldots,m_1$, respectively. 

This analysis for a single fault carries over to multiple faults if they stay spatially separated. Then, in order to apply Corollary \ref{cor1}, we set $I^1=I^1_m=\{0\}$, $p^1=1$, which reflects the presence of the reliable server $\mathcal{S}_0$ in the compute network. Next, we introduce index sets $I^s$, $s=2,\ldots,S-1$, of size $|I^s|=l+1$ with $p^s=l$ by identifying them with the neighborhoods $J$ of the currently failed  $\mathcal{C}_i$, i.e., we have $S-2$ compute nodes still in fail state for which failure occurred in some previous iteration step.  Our assumption of spatial fault separation implies that these $I^s$ are mutually disjoint. Finally, $I^S$ contains the remaining indices, and the random set $I^S\backslash I_m^S$ of size $f'_m$ corresponds to the newly failing $\mathcal{C}_i$ in the current iteration step. Consequently, $f_m=S-2+f_m'$. The number of correctly working $\mathcal{C}_i$ with $i\in I^S$ is 
$$
p^S=p_m-1-(S-2)l=n-f_m-(S-2)l=n-f'_m-(S-2)(l+1)
$$ 
and $|I^S|=n-(S-2)(l+1)$. We thus find that
\be\label{MultF}
\underline{r}_p=\min(\frac{l}{l+1},\frac{n-f'_m-(S-2)(l+1)}{n-(S-2)(l+1)}),\qquad \bar{r}_p=1,
\ee
and the application of Corollary \ref{cor1} provides the corresponding bound for the error reduction in one iteration step. If the
number $S-2$ of failed nodes is moderate and there are almost no newly failing nodes, this bound will be identical to 
(\ref{A33}) since then $\underline{r}_p=l/(l+1)$. In the extreme case of no failed nodes from
previous iteration steps, e.g., when failing nodes can be restarted immediately, we have $S=2$, and $p^2=p_m-1$ is the number of correctly working $\mathcal{C}_i$, $i=1,\ldots,n$. In this situation, the error reduction factor in the bound will be $(1-(p_m-1)^2/(n^2\kappa))$ if the appropriate $\xi_m$ is used. All these bounds
will also hold if $\xi_m$ is selected according to
the steepest descent rule (\ref{XiOpt}). Although the application of
Corollary \ref{cor1} only gives a crude upper bound for the expected error decay per step, it shows that both, the storage redundancy characterized by the integer $l\ge 1$ and the failure rate $f_m/n$ of faulty $\mathcal{C}_i$, may
impact the actual convergence behavior.

The following numerical experiment highlights the influence of the amount of redundant storage at the $\mathcal{C}_i$ characterized by the parameter $l$.
For each $i=1,\ldots,n$ we define beforehand a neighborhood of $\Omega_i$ consisting of $l$ neighbors, where $l=1,2,...,8$ is fixed (for $d=2$ larger values $l>\bar{l}=8$ are prohibitive in practice as they would increase the communication cost considerably).
For each node $\mathcal{C}_i$  we repeatedly run an independent Weibull process (for simplicity with the same parameters $k_1$, $\lambda_1$ for each node) that indicates the begin of 
a failure, followed by another Weibull process for the length of the failure, again independently and with the same parameters $k_2$, $\lambda_2$
for each node. The parameters of these two processes will be tuned such that in each compute cycle a small percentage of the nodes is failing (below we will use realizations with about $10\%$ and $1.5\%$ average per cycle failure rate, respectively). To achieve this
from the start of the iteration, we will initialize the fault process accordingly (and will not start with all nodes in good condition).
During any cycle, for each of the nodes in failed state, the random index pick from the associated index set $J$ given in (\ref{Jset}) is implemented, this 
gives the set $I_m$ for this cycle. Since we cannot guarantee spatial separation of failing nodes in this scheme, there might be conflicts
which are ignored (this is, once a node is asked to switch to a failing node's subproblem solve it will do so and ignore later requests of 
other failing nodes). This is a worst case scenario that can be implemented without additional communication overhead. Recall that the coarse subproblem is solved in each iteration cycle by the reliable server $\mathcal{S}_0$.

\begin{figure}[ht]
\centering
\includegraphics[width=0.46\textwidth,height=0.26\textheight]{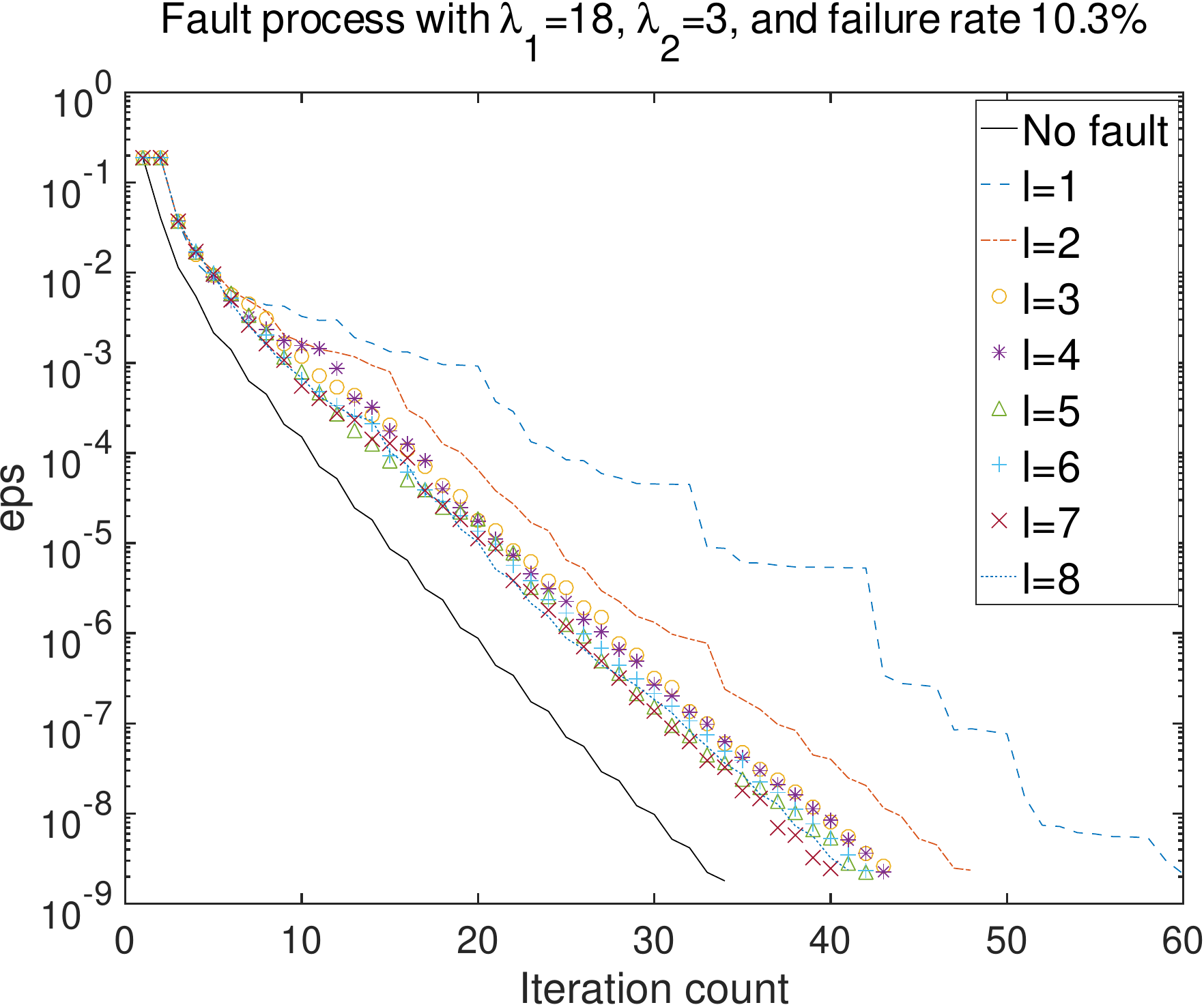} $~~~~~~~~~~~$
\includegraphics[width=0.46\textwidth,height=0.26\textheight]{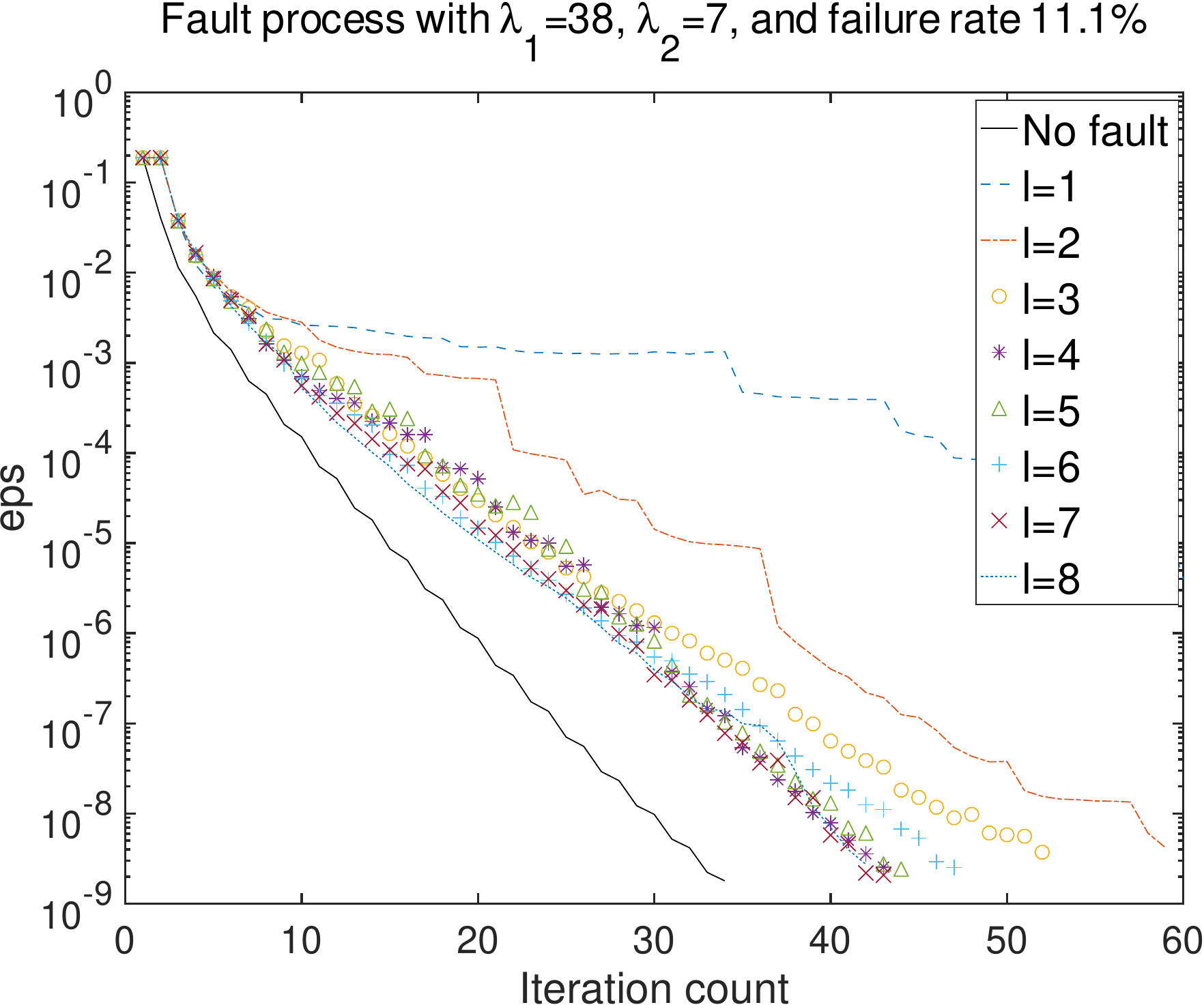}
\caption{Influence of the size  $l$ of subdomain neighbor groups $J$ for fault scenarios generated by two different Weibull processes  with an average per cycle failure rate of about $10\%$ on the convergence of the method (\ref{Rec0})
with $\xi_m$ chosen by the steepest descent rule (\ref{XiOpt}). For comparison, the solid line shows the error behavior of the method in a no-fault situation.}
\label{Fig5}
\end{figure}
According to our theoretical estimate (\ref{A33}), we expect a certain deterioration of convergence speed for small 
$l$, especially if the failure of a few nodes can last for many cycles. To study this aspect, we took again the test problem
from subsection \ref{sec13}, i.e., we set $a(x)=f(x)=1$ in (\ref{MP}) and use an overlapping DD space splitting with parameters $n_1=400$, $n_0=20$, $\ell=6$. We employ the steepest descent rule (\ref{XiOpt}) to determine $\xi_m$. In order to be able to watch the iteration for more cycles and to check if the Weibull parameters used to simulate the fault process impact the
convergence behavior, we chose the smaller value $\epsilon_0=10^{-8}$ as termination criterion. The standard additive Schwarz method took $34$ iterations to reach this relative error reduction. 

\begin{figure}[ht]
\centering
\includegraphics[width=0.46\textwidth,height=0.26\textheight]{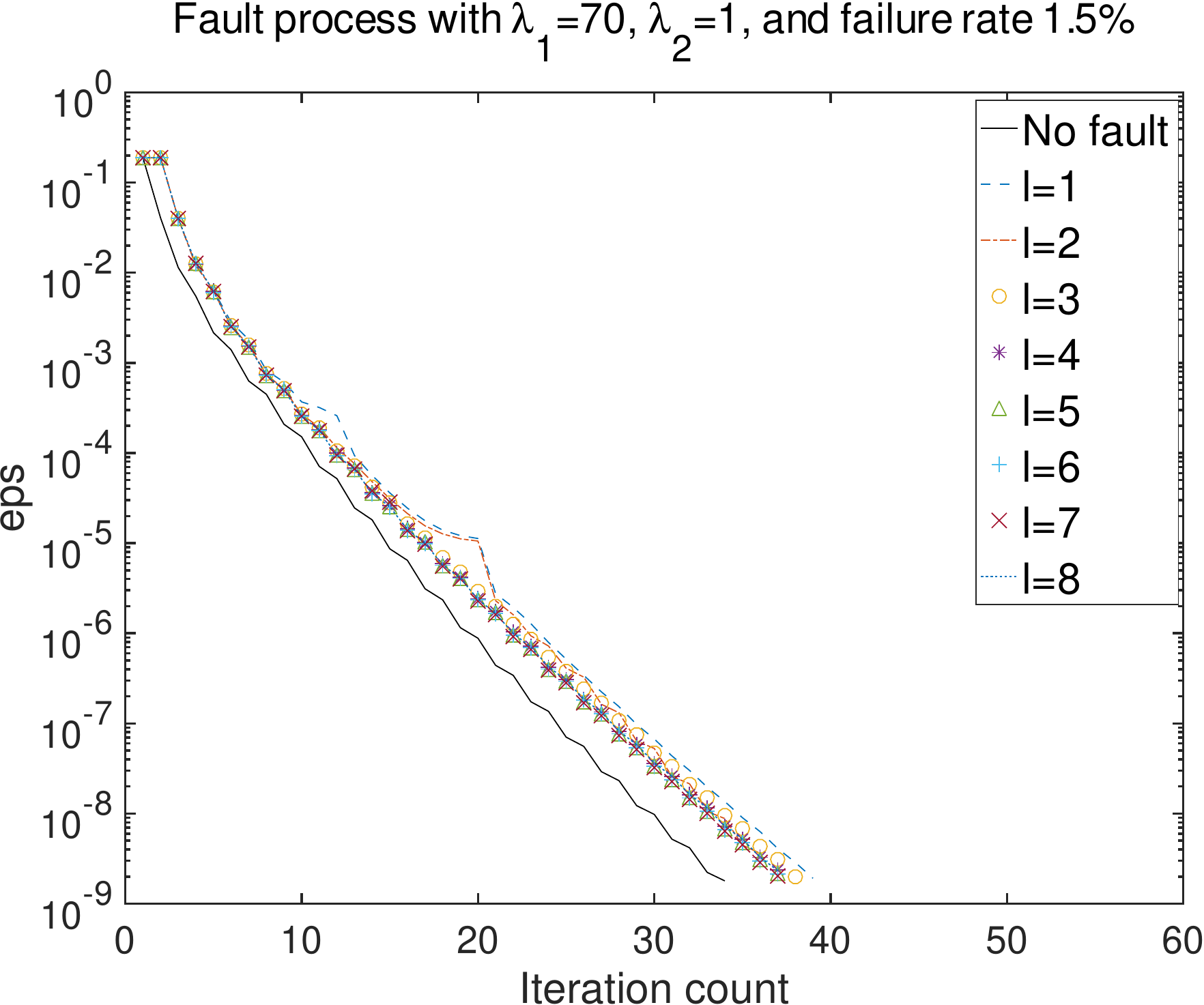} $~~~~~~~~~~~$
\includegraphics[width=0.46\textwidth,height=0.26\textheight]{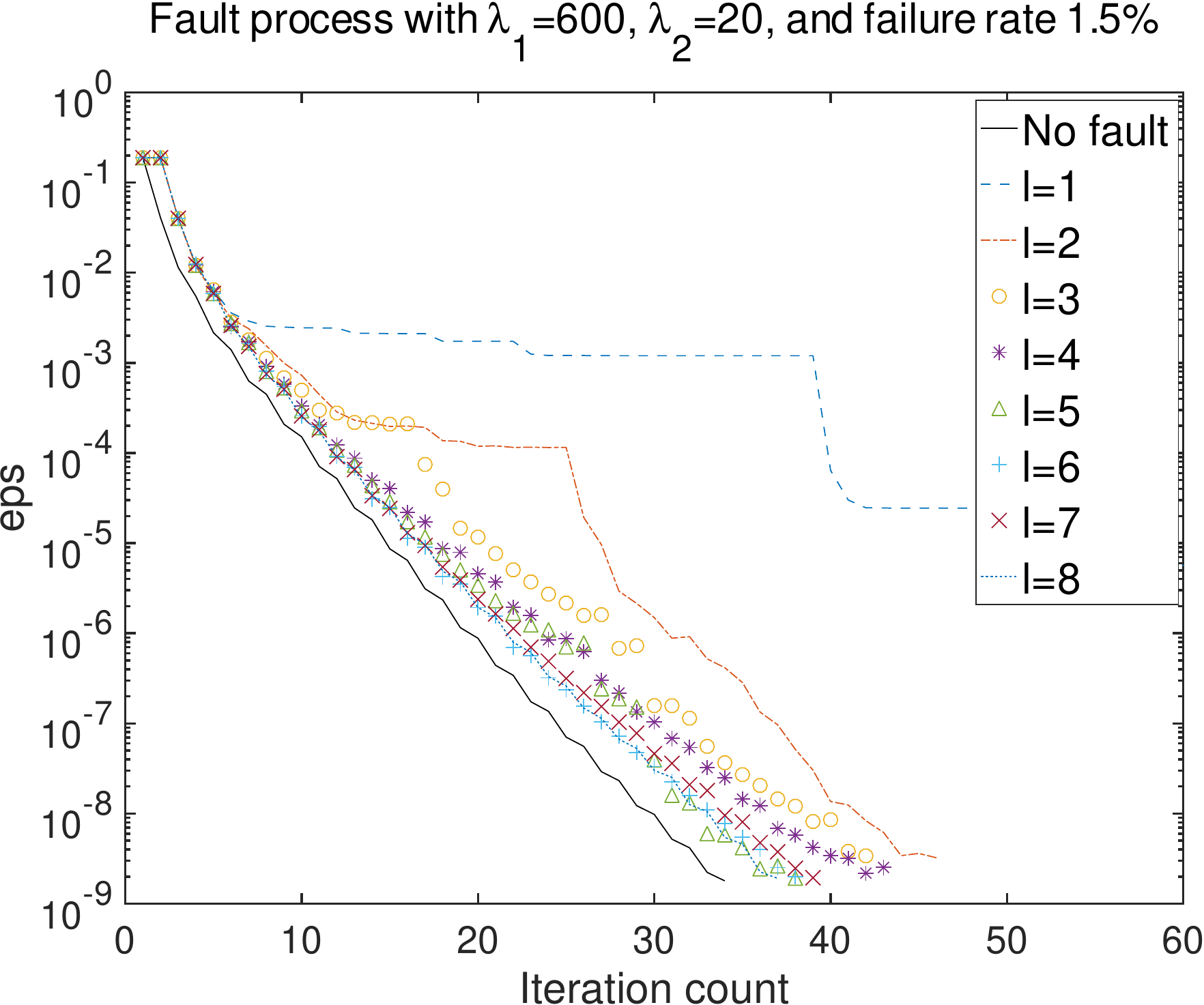}
\caption{Influence of the size  $l$ of subdomain neighbor groups $J$ for fault scenarios generated by two different Weibull processes with an average per cycle failure rate of about $1.5\%$ on the convergence of the method (\ref{Rec0})
with $\xi_m$ chosen by the steepest descent rule (\ref{XiOpt}). For comparison, the solid line shows the error behavior of the method in a no-fault situation.}
\label{Fig6}
\end{figure}

Figures \ref{Fig5} and \ref{Fig6} show the results of runs with different fault processes. In all four cases, the Weibull shape parameter for the time to the next failure of a compute node was set to $k_1=0.5$ (such a value was also used in \cite{PAS2014}) while for the length of the failure the value was set to $k_2=1$ (exponential distribution). The associated scale parameters were chosen as follows: For the experiment in 
Figure \ref{Fig5} (left), we took $\lambda_1=18$ for the time to the next failure (meaning that on average a compute node fails after $36$ cycles) and $\lambda_2=3$ for the length of failure (meaning that on average a node remains failing for $3$ cycles only). These choices resulted in a realization where in each compute cycle on average about $10.3\%$ of the $n=400$ compute nodes 
were failing. For the experiment in Figure \ref{Fig5} (right)  the corresponding values were $\lambda_1=38$, $\lambda_2=7$ resulting in a realization with an average per cycle failure rate of about $11.1\%$. In the latter case, an individual node remains non-faulty on average for $76$ cycles but when it fails it stays faulty on average for $7$ cycles. This is the situation in which we expect to see a more profound impact of the parameter $l$. We here deliberately considered
relatively large average per cycle failure rates to show that our approach is robust and behaves according to the theoretical predictions.

The results for different $l=1,\ldots,8$ depicted 
in the graphs are in complete agreement with our theory, in particular, they confirm the estimate (\ref{A33}) qualitatively: Larger values for the size $l$ of subdomain neighbor groups lead to better performance. This can also be seen from the number of iterations for each of the failure scenarios given in Table \ref{Tab2}: As expected, 
the influence of $l$ becomes more visible with the increase of $\lambda_2$, i.e., with the number of cycles a failed compute node
remains in failing state.  The more irregular convergence behavior in the graphs for $l=1,2$, especially for larger values of $\lambda_2$, may have different reasons.
On the one hand, it may be due to a certain loss of spatial separability of faulty compute nodes. Our above mentioned crude conflict resolution strategy may have resulted in neglecting certain subproblems for many cycles. On the other hand, the graphs depict computable error indicators as explained in the appendix, not the exact errors $\|e^{(m)}\|$, which may also result in a less smooth error decay.

In Figure \ref{Fig6}, we show test results for fault processes with smaller 
average per cycle failure rates of about $1.5\%$. Such a failure rate is more
realistic given the current predictions for failure rates of processors in large parallel architectures. For Figure \ref{Fig6} (left) the Weibull parameters were
$\lambda_1=70$, $\lambda_2=1$, which means that, on average, individual compute nodes stay alive for $140$ cycles but are restarted almost immediately after a failure. 
The parameters for Figure \ref{Fig6} (right) were $\lambda_1=600$, $\lambda_2=20$,
which represents the other extreme: Relatively long failure times of compute nodes
after very long periods of correct functioning. The results again  confirm the 
predicted dependence of convergence rates on the parameter $l$ characterizing the 
amount of redundancy which increases with the parameter $\lambda_2$ characterizing the average failure time of compute nodes. Compared to Figure \ref{Fig5}, we also see
the impact of the average per cycle failure rate $r_f$: The smaller $r_f$, the closer the performance gets to that of the additive Schwarz iteration in a non-faulty compute network (at least, if $l$ is large enough). This is also illustrated by the iteration counts for the above four test cases recorded in Table \ref{Tab2}.
\begin{table}[ht]
\centering
\begin{tabular}{cc|rrrrrrrr}
&\multicolumn{8}{c}{Iteration counts for different $l$}\\
\hline
$\lambda_1$& $\lambda_2$ & $l=1$ & $l=2$ & $l=3$& $l=4$& $l=5$& $l=6$& $l=7$& $l=8$\\
\hline
18 & 3 & $60\;$ &$48\;$ &$43\;$ & $43\;$&
$42\;$ &$42\;$ &$40\;$ & $41\;$\\
38 & 7 & $>100\;$ &$59\;$ &$52\;$ & $43\;$&
$44\;$ &$47\;$ &$43\;$ & $42\;$\\
70 & 1 & $39\;$ &$36\;$ &$38\;$ & $37\;$&
$37\;$ &$37\;$ &$37\;$ & $37\;$\\
600 & 20 & $91\;$ &$46\;$ &$42\;$ & $43\;$&
$38\;$ &$38\;$ &$39\;$ & $37\;$
\end{tabular}
\caption{Iteration counts for reaching a relative error reduction of $\epsilon_0=10^{-8}$ for the fault scenarios used for the graphs in Figure \ref{Fig5} and \ref{Fig6}. The corresponding Weibull scale parameters $\lambda_1$, $\lambda_2$ are shown in the table, the Weibull shape parameters are $k_1=0.5$, $k_2=1$.
The iteration with no faults needed $34$ iteration steps to termination.}
\label{Tab2}
\end{table}

To summarize, even under the assumption of a compute network with predominantly local communication and distributed data storage, we can still
get reasonable convergence rates if we allow for redundant storage at the unreliable compute nodes with slightly larger values $l$ than proposed in \cite{CXZ2017}, and treat the coarse problem
at a reliable server. We refer to the appendix for some hints on implementation details.

As a final remark, let us mention that similar considerations are possible for other compute network architectures.
For instance, in \cite{RiMo2017} a server-client architecture was used to
achieve fault tolerance for an overlapping DD method for (\ref{MP}) in one and two dimensions 
without global communication. To this end, the PDE problem is turned  into a fixed-point formulation for the 
system of local boundary-to-boundary maps for the restrictions of the solutions to the subdomain boundaries
$\partial \Omega_i$ and their neighbors $\partial \Omega_j \cap \Omega_i$.
For those, approximations are generated assuming a fault model, where faults are not lost subproblem solves but may be accidentally missing data for the boundary-to-boundary maps. In essence, this represents a reformulation of an overlapping DD splitting, similar to (\ref{DDS}) but without a coarse space $V_0$, whose convergence properties will obviously deteriorate with the number of subproblems. The boundary-to-boundary maps
are executed on unreliable clients attached to a network of reliable servers. We refer to \cite{RiMo2017,SaRi2017} for details. 

We use this server-client model of \cite{RiMo2017} to discuss one more potential application of our results in section \ref{sec1}. Think of a network of reliable servers $\mathcal{S}_0, \mathcal{S}_1,\ldots, \mathcal{S}_L$, and consider,  besides the coarse partition $T_0$, another, still coarser overlapping partition of $\Omega$ into $L<<n$ domains
$\tilde{\Omega_j}$, each of which is the union of about $n/L$ subdomains $\Omega_i$ (each $\Omega_i$ belongs to exactly one $\tilde{\Omega}_j$). Each of the $L$ servers $\mathcal{S}_j$, $j=1,\ldots,L$,  has enough compute power and memory to keep safe copies of the distributed static and dynamic data arrays
associated with the subdomains $\Omega_i$ forming $\tilde{\Omega}_j$. This way each $V_i$ subproblem, $i=1,\ldots,n$, is owned by exactly one server.  Each of these $L$ servers has clients with point-to-point communication to the server
but not with each other that will deal with solving the subproblems owned by the server (in other words, each server with its clients represents a master-slave architecture, similar to subsection \ref{sec13}). We assume that the number of clients per server is such that during a compute cycle each subproblem owned by a server can be assigned to one of its clients, and that there is no 
correlation between client failure and subproblem assignment. As in the case of the local communication network, the server $\mathcal{S}_0$ is reserved for
dealing with the $V_0$ subproblem solve and global error computations.
Communication between the servers is considered reliable, $\mathcal{S}_0$ needs to communicate with all other $\mathcal{S}_j$ while the servers $\mathcal{S}_j$, $j=1,\ldots,L$,
need to be linked with $\mathcal{S}_0$ and in addition with those $\mathcal{S}_{j'}$ for which 
$\tilde{\Omega}_j\cap \tilde{\Omega}_{j'}\neq \emptyset$. We call such $\mathcal{S}_j$ and
$\mathcal{S}_{j'}$ neighboring. Note again at this point that the overall communication of necessary data between servers $\mathcal{S}_j$, $j=1,\ldots,L$, can be performed in parallel in $\approx$ $\bar{\ell}'$ sweeps, where $\bar{\ell}'$ denotes now an upper bound for the number of neighboring servers of any given $\mathcal{S}_j$. 

Under these assumptions, we can use the same strategy as in subsection \ref{sec13} locally
for each server-client subgroup of the network. The only difference is that, in addition to feeding the clients with randomly assigned
static and dynamic data arrays for subproblem solves, each server communicates dynamic data for solving the  $V_0$ problem and computing error estimators to the server $\mathcal{S}_0$. 
This is a small difference compared to subsection \ref{sec13}, where we silently assumed that all subproblem solves including the coarse subproblem are of the same run-time complexity and therefore could be dealt with by the slave nodes, even though one could have delegated the subproblem solve associated with $V_0$ to the master node as well. Note that the point-to-point communication
between the servers $\mathcal{S}_j$, $j=1,\ldots,L$, and the specialized server $\mathcal{S}_0$ involves only arrays of small size proportional to the number of
subdomains assigned to each server. We refer to the appendix for some more details.


\subsection{Estimate of parallelization gains}\label{sec15} 
We now give a rough calculation of the runtime per cycle for the three architectures discussed in the previous subsections based on
a simplified runtime model for the involved subproblem solves, update, and communication steps. Note that a typical cycle consists of a communication step between reliable server(s) and compute nodes or between compute nodes, a parallel subproblem solve step, followed by a reverse communication step
and the update step.

For the subproblem solve step, we assume that one subproblem 
solve (provided that dynamic and static data arrays associated with $\Omega_i$ are made available) on one compute node
takes $C_{s}M_i$ units of cpu-time, i.e., it scales linearly in the dimension $M_i\approx M$ of the $V_i$ subproblem with a constant $C_s$ that depends
on the solver, the required accuracy, and the cpu-speed, and may be large. For PDE problems such as our model problem (\ref{MP}) such an assumption is realistic if optimized multigrid methods are applicable to the subproblems. If suboptimal local solvers (e.g., direct solvers for sparse linear algebra problems) are used instead, this may not be true, in which case some of the conclusions below will look more optimistic since it is essentially only the solve step that can be fully parallelized. Similarly, the update step which consists of performing linear combinations of vectors stored in distributed format for
each $\Omega_i$, merging information received from neighboring $\Omega_{i'}$, and computing local contributions to error estimators takes
$C_u M$ units of cpu-time, where typically $C_u<<C_s$ can be assumed.  These assumptions will also be applied to the cpu-time of the coarse problem solve with $M_i\approx M$ replaced by $M_0\approx n$.

As to the communication steps, the cpu-time model is $C_{0c}+C_c s$ for a one-to-one communication of an array of size $s$
between two compute nodes or servers, where $C_{0c}$ is the absolute time for opening the connection, and $C_c s$ is the transmission time for the actual data. Scheduling of many connections from or to any given node is sequential.
With these model assumptions at hand, we will derive rough runtime estimates for the parallel DD algorithms on the three architectures. 
Recall that the size of the whole problem, i.e., the dimension of $V$, is $N=\dim V\approx n M$. Given the theoretical option of solving the problem with the same cost model on a single, large and reliable computer, we could achieve the solution of the problem in time $(C_s+ C_u)nM$ (this does not include the precomputation of $A$ and $b$ and other static information).

We start with the \emph{master-slave network} of subsection \ref{sec13}. Even though we assumed any value $p\le n+1$ of available slaves, to be somewhat compatible with the other scenarios, we look at the case $p= n+1$. In each cycle, we send in an one-to-all communication step static and dynamic data arrays
associated with all subproblems including the $V_0$ subproblem in a random assignment from the master to the $n+1$ slaves in total time $(n+1)C_{0c}+C_c(n M+M_0)$,
then solve the subproblems in time $C_s\max( M,M_0)$ at the slaves, and communicate the results back in an all-to-one communication step to the master node,
again in time $(n+1)C_{0c}+C_c(n M+M_0)$, whereas constants may be different. At the master, the distributed data representation needs to be synchronized at cost $C_u (nM+M_0)$. Using the fact that $M_0\approx n$, this results in a time budget of approximately
\be\label{TMS}
T=2(n+1)C_{0c}+C_{c}n(M+1) + C_s\max( M,n) +C_u n(M+1)
\ee
per cycle. In comparison, if the whole iteration step would have been performed by the master node, we would obtain a time budget of $(C_s+ C_u)Mn$. Thus, other than having freed the master node from the serial solve step, there is no gain from outsourcing computational work to the unreliable slaves unless $C_s$ is much larger than the constants associated with update and communication steps.

For the \emph{local communication network} described in subsection \ref{sec14}, a possible implementation is described in the appendix,
we refer to it. The solve step 1 takes $C_s \max(M,M_0)$, and the update step 3 takes $C_u M$ units of time, respectively (the constant $C_u$ adsorbs the slightly increased amount of work of the neighbors of a failing node). In the communication steps 2 and 4, we have a local part which costs approximately $\bar{l}(C_{0c}+C_c M)$
units of time, where $\bar{l}$ is given by (\ref{Locall}),
and we have the one-to-all and all-to-one communication of small amounts
of data to the specialized compute node responsible for the coarse subproblem. Since $M_0\approx n$, the latter takes $n(C_{0c}+C_c)$ units of time in the worst case
(again, constants may be different).
Altogether, this results in an overall time budget of
\be\label{TLC}
T=2\bar{l}(C_{0c}+C_{c}M)+ 2n(C_{0c}+C_{c}) + C_s \max(M,n) + C_u M, 
\ee
per cycle. This estimate is linear in $M$ and $n$, possibly with network- and cpu-dependent constants
for solve, update, and communication steps but may grow with the amount of local overlap in the underlying 
domain partition $\{\Omega_i\}$. It also depends on the dimension-dependent constant $\bar{l}$.  From the point of view of parallel efficiency, this is as good as one can expect if a coarse subproblem is included.

We note that the possible gain of the approach in subsection \ref{sec14} is visible only if $C_s$ is large compared to other constants, i.e., if the solve time in step 1 is dominating. Otherwise, we could modify the algorithm outlined in the appendix as follows: In the case of a failing node  $\mathcal{C}_i$, any of the $l$ designated neighbors $\mathcal{C}_{i'}$ with current copies of data arrays for the $V_i$ subproblem solve takes on the responsibility for solving the $V_i$ subproblem in addition to its own subproblem. This results in a solve time of roughly $2C_s M$ in each cycle while leaving all other steps unchanged. If the solve time is indeed dominating, this alternative ABFT
approach would therefore double the runtime while the expected slowdown due to the ABFT approach described in subsection \ref{sec14}
depends on the amount of local storage redundancy characterized by $l$, and becomes less visible with larger values of $l$. The same tradeoff has been discussed for $l=1$ in a similar situation
in \cite{CXZ2017}.

Finally, for the \emph{server-client network} mentioned at the end of subsection \ref{sec14}, the overall time is intermediate to the previous two cases. The difference to the master-slave network of subsection \ref{sec13} is that the role of the master computer is now played by a reliable network consisting of a special server $\mathcal{S}_0$, and of $L$ servers $\mathcal{S}_j$, $j=1,\ldots,L$, each of which keeps the data related to about $n/L$ subproblems.
During a cycle, each server $\mathcal{S}_j$ ($j\neq 0$) acts as a small master and feeds its $\approx n/L$ slaves, taking $2(n/L)(C_{0c}+C_c M)$
units of time for serial communication to and from its clients, and $C_s M$ units of time for the parallel subproblem solve step. During the subproblem solve step at the clients, the servers $\mathcal{S}_j$, $j=1,\ldots,L$,
can do the communication with the special server $\mathcal{S}_0$ responsible for the coarse problem associated with $V_0$, the computation of error indicators, and the maintenance of information about the number $p_m$ of correctly functioning clients (this involves an all-to-one and one-to-all communication with $\mathcal{S}_0$ costing $2L(C_{0c}+C_c M_0/L)$ units of time and a solve step of cost $C_s\max(M,M_0)$). If the communication in the server network is sufficiently fast compared to the solve time at the clients, this time may not matter. After the $L$ servers 
have received data of the subproblem solves from all correctly working clients (and coarse problem and global data from the special server $\mathcal{S}_0$),
they independently perform an update step costing $C_u (nM+M_0)/L$ units of time. This is followed by a synchronization step between neighboring servers
(this is the transmission of dynamic data arrays such as $x_{ii'}$ associated with the overlap region of subdomains $\Omega_i$ and $\Omega_{i'}$ owned by different servers). For large $L$
and again with $M_0\approx n$,  this cost is negligible compared to the already accrued overall cost of
\be\label{TSC}
T=2(n/L)(C_{0c}+C_c M)+ 2L(C_{0c}+C_c n/L) + C_s \max(M,n) + C_u (M+1)n/L,
\ee
where we can safely  assume that $L<< n$. Compared to the budget (\ref{TMS}) for the master-slave network, we benefit by reducing communication and update time by roughly the factor $L$, i.e., the size of the server network. To come close to the estimate (\ref{TLC}),
one would need to assume small values $n/L$. This situation is, however, not of practical interest because it implies the presence of a large reliable server network with only few compute nodes attached to each of the servers.

\section{Concluding remarks and future work}\label{sec4}
In this article, we have considered an example for algorithm-based fault tolerance (ABFT), namely how to make domain decomposition methods in PDE
applications more fault tolerant. To this end, we considered stochastic subspace correction algorithms
and developed a general theoretical foundation for their convergence rates under weak randomness and independence assumptions for failure of subproblem solves. As an application, we used a standard overlapping domain decomposition
method for a simple two-dimensional Poisson problem and showed that our convergence theory for stochastic subspace
correction methods indeed gives proven convergence rates also in the faulty case and results in
the design of fault-tolerant methods, e.g. for local communication networks, with quasi-optimal parallel cost complexities.

So far, we employed our theory to a simple two-dimensional model problem only and not to a large three-dimensional, time-dependent, nonlinear real-life simulation problem from e.g. physics or engineering yet. 
This is future work. We also used in our experiments fault rates up to $15$ percent which is unrealistic and much too high. 
Nevertheless, this demonstrated the robustness of our proven convergence bounds also in such a situation and it is clear 
that in more practical situations, i.e., for smaller failure rates in the per-mil range, the resulting convergence rates must be nearly as good as the ones in the non-faulty case.  

Note at this point that our theory can be applied to other space splittings as well and will then lead to associated convergence bounds and corresponding fault-tolerant parallel algorithms. One example is the case of inexact subproblem solvers, another one are various multilevel and multigrid solvers which might be analyzed as subspace correction methods in a similar way. 
A further example may be the sparse grid discretization \cite{BG04} and the so-called combination approach for higher-dimensional partial differential equations, e.g. in its original version \cite{combi} or in its improved version 
as the so-called Opticom \cite{gaheg}, see also the results in \cite{GHO15}. It involves a combination of smaller, in general non-isotropic discretizations of the problem at hand which can be treated completely independent of each other, each e.g. by a parallel DD method itself. This way a second level of parallelization is introduced which altogether leads for high-dimensional problems to a huge amount of decoupled subproblems. A main issue is again if a certain amount of subproblem solvers is faulty.   Then, our theory also gives results for fault-tolerant versions of such algorithms in a straightforward way. An example for a real life application under study involve the gyrokinetic equations for microturbulence to compute gyroradius-scale fluctuations and the resulting transport coefficients in magnetized fusion/astrophysical plasma, see e.g. \cite{HHOBP17,HHHB16} and the references cited therein.

Note at last that another promising area of application for our theoretical results on stochastic subspace correction algorithms is machine learning and data analysis. In such applications, there are many situations where randomness in the subproblem selection is just given by the problem and the data under consideration. 
An example would be online learning algorithms for solving classification or regression problems in a reproducing kernel Hilbert space setting which again can be interpreted in the framework of Schwarz methods. There, the subproblems can be associated with the samples in the training set which makes the incremental learning algorithm automatically randomized. Also in such situations our theory can be applied in a straightforward manner. 

\section*{Acknowledgment}
M. Griebel acknowledges the support from the DFG priority program 1648 "Software for Exascale Computing" within the project "EXAHD - An Exa-Scalable Two-Level Sparse Grid Approach for Higher-Dimensional Problems in Plasma Physics and Beyond".
The main results of this paper were obtained during a yearlong stay of P. Oswald at the Institute for Numerical Simulation (INS) sponsored by the
Hausdorff Center for Mathematics of the University of Bonn and funded by the Deutsche Forschungsgemeinschaft. He is grateful for this support.

\bibliographystyle{amsplain}

\section{Appendix: Distributed implementation} \label{sec3}

With the space splitting (\ref{DDS}) fixed and an implementation on a local communication network as described in subsection \ref{sec14} in mind,
we introduce the following notation for static and dynamic data arrays used below. Denote by $x$ the nodal basis vector of size $N=\dim(V)$ representing 
an arbitrary $u\in V$. For any given $i$, $x_i$ denotes the subvector of $x$ of size $M_i$ corresponding to the nodal basis in $V_i$ which is a subset of the nodal basis in $V$. The associated index set is denoted as $J_i$. We also need the subvectors $x_{ii'}$ of $x_i$ of length $M_{ii'}<M_i$
which correspond to the basis functions  in  $V_{ii'}:=V_i\cap V_{i'}$ (obviously, for each $i$ only a small number of nontrivial $x_{ii'}$ need to be considered). Note that, even though formally $x_{ii'}=x_{i'i}$,  in the actual iteration the vectors ${x}_{ii'}$ and ${x}_{i'i}$ may differ temporarily.
E.g., if $x_i$ is updated in a $V_i$ subproblem solve at $\mathcal{C}_i$ then $x_{ii'}$ changes while the neighbors $\mathcal{C}_{i'}$ may have different
$x_{i'i}$ from the previous iteration or their own subproblem solve step. The 
equality $x_{ii'}=x_{i'i}$ is again guaranteed after the next communication step and the update steps are performed.
Obviously, $x$ can always be reconstructed from its distributed representations
$\{{x}_i\}$.  The extension-by-zero maps $x_i\in \mathbb{R}^{M_i} \to x=(x_i,0)\in \mathbb{R}^{N}$ define the operators $R_i$ used for the space splitting (\ref{DDS}), $i = 1,\ldots,n$. Below, we will use the same notation for the operator $R_i$ and its matrix representation. The restriction map $R_i^T:\,x\in \mathbb{R}^{N} \to x_i\in \mathbb{R}^{M_i}$ defines the adjoint to $R_i$, $i = 1,\ldots,n$.
It is convenient to precompute some sparse matrices and vectors. Let $A$ be the nodal basis discretization matrix,   
and let $b$ be the right-hand side of the sparse linear system representing the FE discretization of (\ref{MP}) associated with $V$. Then, for $i=1,\ldots,n$, we introduce with $A_i=R_i^TAR_i$ the overlapping diagonal blocks of $A$ of size $M_i\times M_i$ associated with $V_i$. Similarly,
we have $b_i:=R_i^Tb$, while $A_{ii'}=R^T_{i'}AR_i$ are the submatrices of size $M_{i'}\times M_i$ associated with the overlap regions
$\Omega_i\cap \Omega_{i'}$.
Furthermore, denote by $A_0$ and $b_0$ the stiffness matrix and right-hand side of the nodal basis discretization associated with $V_0$,
respectively. Due to our assumptions, their representation can be produced from $A$ and $b$ by the formulas
$A_0=R_0^TAR_0$ and $b=R^T_0b$, where $R_0: V_0 \to V$ is the natural embedding operator whose distributed matrix representation $\{R_{0i}\}_{i=1,\ldots,n}$
is given as follows: Let $x_0$ denote the nodal basis vector of an element $u_0\in V_0$. Then, ${R}_{0i}$ is a submatrix of the $N\times M_0$ matrix $R_0$ which corresponds to mapping the subvector $x_{0i}$ of $x_0$  associated with the nodal basis functions in $V_0$ whose support intersects with ${\Omega}_i$ to the subvector $(R_0x_0)_i$ of $R_0x_0$.  Note that in the DD setting the vectors $x_{0i}$ have relatively small and uniformly bounded size, independently of $M$ and $n$, while the size of $x_0$ equals the dimension $M_0$ of $V_0$, and scales linearly with $n$. 
Finally, to run the one-step recursion (\ref{Rec0}) and the
accelerated method (\ref{u}-\ref{v}), a certain number of small arrays containing parameters and auxiliary data are needed. By $\xi$ we will denote a vector of length
$n+1$ whose entries contain the relaxation parameters
$\xi_i:=\xi\omega_i$, $i=0,1,\ldots,n$, while $e'$ denotes a vector for storing local error indicators associated with the subproblems.

Using this notation, 
we next give some implementation details for the iteration (\ref{Rec0}). The changes for the implementation of the accelerated iteration (\ref{u}-\ref{v}) are briefly mentioned at the end. We concentrate on the local communication computer network of subsection \ref{sec14}
which seems most promising from the point of view of overall efficiency.
In the following, we will silently add vectors and perform matrix-vector products with different index sets 
by padding the vectors with zeros to the correct dimensions.

Let us first write the iteration in vector matrix notation as a sequence of single update steps with one subproblem solve at a time. We denote by $x$ and $x_{new}$ the nodal basis vector of the iterate $u$ before and after the update, respectively.
Similar notation is used for the residual $r=b-Ax$, its coarse grid projection $r_0=R_0^Tr$, and the distributed representations thereof.
If the subproblem has index $i\neq 0$ then
$$
x_{new} = x + \xi_i R_i A_i^{-1}R_i^T(b-Ax) = x+ \xi_i A_i^{-1}r_i,\qquad r_i=R_i^Tr,
$$
and thus in distributed format
\bea
x_{new,i'}&=& x_{i'} + \xi_i R_{ii'}d_i,  \qquad R_{ii'}:=R^T_{i'}R_i,\quad d_i:=A_i^{-1}r_i, \\
r_{new,i'} &=& r_{i'} - \xi_i R_{i'}^TAR_i d_i = r_{i'}-\xi_i A_{ii'}d_i,
\eea
where $i'=1,\ldots,n$ (since $x$ can be recovered from the $x_i$ with $i=1,\ldots,n$, updates of $x_0$ are not needed). Recall that $R_{i'}^TR_i$ 
just represents the restriction of $x_i$ to $x_{ii'}$, that the submatrices $A_{ii'}=R_{i'}^TAR_i$ of $A$ are non-zero for only a few $i'$,
and that the vector $A_{i0}d_i$ has uniformly bounded size.
Thus, possibly except for the solution of the problem $A_jd_j=r_j$, computation time and storage for 
the relevant quantities $x_{ii'}:= \xi_i R_{ii'}d_i$ and $r_{ii'}:= \xi_i A_{ii'}d_i$ will remain roughly proportional to $M_i\approx M$. Note that in case $i'=i$
the update is
$$
r_{new,i}= r_{i}-\xi_i A_{ii}d_i=r_i-\xi_i A_i A_i^{-1}r_i=(1-\xi_i)r_i,
$$
which could be used to check the validity of returned results. In the implementation, the updates for $i'=i$ ($x_{new,i}$
and $r_{new,i}$), and $i'=0$ ($r_{new,0}$) will be performed during the solve step at the compute node responsible for  this subproblem. The update to $x_{new,i'}$ and $r_{new,i'}$ with $i'\neq i$  requires communication with the neighboring compute nodes.

If $i=0$ (i.e., the subproblem is associated with the coarse space $V_0$, and is executed on the reliable server $\mathcal{S}_0$) we have for $i'=0$
$$
r_{new,0} = r_{0} - \xi_0 R_{0}^TAR_0 A^{-1}_0r_0=(1-\xi_0)r_0 ,
$$
while for  $i'=1,\ldots,n$ the update formula is 
\bea
x_{new,i'}&=& x_{i'} + \xi_0 R_{0i'}d_0,  \qquad R_{0i'}:=R^T_{i'}R_0,\quad d_0:=A_0^{-1}r_0, \\
r_{new,i'} &=& r_{i'} - \xi_0 R_{i'}^TAR_0 d_0 = r_{i'}-\xi_0 A_{0i'}d_0.
\eea
Note that $R_{0i'}=R_{i'0}^T$ and $A_{0i'}=A_{i'0}^T$ need only the small subvector $d_{0i'}$ associated with the few nodal basis functions with support intersecting $\Omega_{i'}$.
I.e., once $d_0=A_0^{-1}r_0$ is
computed, each evaluation of $x_{0i'}:= \xi_0 R_{0i'}d_0$ and $r_{0i'}:= \xi_0 A_{ii'}d_0$ is of complexity proportional to $M_{i'}\approx M$,
and should be done at the compute node responsible for the subdomain $\Omega_{i'}$ after the associated data $d_{0i'}$ is received from $\mathcal{S}_0$. 

We give now a more detailed description of one cycle for executing the recursion step (\ref{Rec0}) in the case of a DD type space splitting
as described in subsection \ref{sec12}. We assume that the
compute network consists of $n$ compute nodes $\mathcal{C}_i$ and a reliable server $\mathcal{S}_0$. 
All global vectors and matrices such as $A$, $R_0$, $b$, $x$ (placeholder for the nodal basis coefficients of the iterates),
$r=b-Ax$ (the corresponding residual) are represented in distributed format, and stored redundantly at the compute nodes $\mathcal{C}_i$, together
with the network neighbor structure.  
 We assume that the compute node $\mathcal{C}_i$ is responsible for $\Omega_i$, $i=1,\ldots,n$, and stores copies of static and dynamic data arrays associated with $\Omega_i$.
The data associated with $\Omega_i$ is also redundantly stored at $l$ neighboring compute nodes $\mathcal{C}_{j_1},\ldots,\mathcal{C}_{j_l}$ (the index set $J=\{j_1,\ldots,j_l\}$ is fixed but obviously depends on $i$). Among them we select one beforehand,  call it \emph{partner} of $\mathcal{C}_i$, and denote it by $\mathcal{C}^\ast_i$. If $\mathcal{C}_i$ fails, its partner $\mathcal{C}^\ast_i$ will play a special role for organizing the local fault mitigation process described in subsection \ref{sec14}. For simplicity, we assume that faults are spatially isolated.

The server $\mathcal{S}_0$ only needs $A_0$ and network connectivity information to be able to receive from and 
distribute to all other compute nodes small data arrays associated with the coarse problem solve. Our description of a cycle starts with
the parallel solution of subproblems according to a given assignment of subproblems to compute nodes, and ends with a new assignment (or with the decision to stop the iteration).
The initial assignment is identical with the setup of the compute network.  

\bigskip
{\bf Recursion}.
\bi
\item[1] {\it Solve step}. For a $V_i$ problem with $i\neq 0$ assigned to one of the compute nodes, this step includes the
computation of $x_{new,i}$ and $r_{new,i}$ according to the above formulas, and of the data arrays $x_{ii'}:= \xi_i R_{ii'}d_i$ and
$r_{ii'}:= \xi_i A_{ii'}d_i$ needed for updates
at neighboring compute nodes. For the $V_0$ problem, this also includes the assembly of $r_0=R_0^T r$ from its locally computed parts, and the computation of the small arrays $d_{0i'}$, $i'=1,\ldots,n$, representing $d_0=A_0^{-1}r_0$ locally. 
In both cases, the scalar  product $e_i=r_i^T d_i$ is computed.
Moreover, at $\mathcal{S}_0$ the global error indicator
$$
\epsilon = \left(\sum_{i=0}^n e_i\right)^{1/2}
$$ 
is computed from $e_0$ and the values  $e_i$, $i=1,\ldots,n$, available from the previous cycle (see Step 4 below).
\item[2] {\it After-solve communication step}. Each active node responsible for a $V_i$ problem with $i\neq 0$ sets up communication with all its neighbors and to $\mathcal{S}_0$, and collects the relevant data $x_{i'i}$, $r_{i'i}$ ($i'\neq i,0$), $d_{0i}$ for updating $x_i$ and $r_i$, the error indicator $\epsilon$ and the scalars $e'_i$, $i=1,\ldots,n$, needed for computing the global error indicator in the next cycle. In this step, a newly failing compute node $\mathcal{C}_i$ has to be detected by all its neighbors, and the partner $\mathcal{C}^\ast_i$ of the failing node will become known
to all of them. The partner node is added to the neighborhoods of all neighbors of the associated failed node and vice versa, to allow for direct partner-to-neighbor communication in this and future cycles. The set of neighbors of a failing node acts in a specific way guided by the associated partner node until the failing node has been restarted. Small modifications of the data shipped from and to this local group of compute nodes are clear from the explanations in the remaining steps. 
\item[3] {\it Update step}. Each active compute node updates the dynamic data arrays $x_i$ and $r_i$ assigned to it by incorporating all information received from the neighbors and the contribution from the coarse subproblem.  For $\mathcal{C}_i$ that are not partners of a failing compute node, this is the dynamic data for the $V_i$ subproblem permanently assigned to the node under consideration. The partner of any of the failing nodes will in addition update the data assigned to the failing node. I.e., in the previous communication step any such partner must have collected data from its own neighborhood and from the neighborhood of the failing node. 
\item[4] {\it Pre-solve communication step}. Each active node sets up communication links with its neighborhood (in the case of the partner of a failing node this is the extended  neighborhood containing also all neighbors of the failing node), and sends the updated dynamic data for synchronization. Thus, for partner nodes $\mathcal{C}^\ast_i$ 
of failing nodes $\mathcal{C}_i$, the amount of data shipped may be temporarily slightly larger. \\
Next, the (non-faulty) compute nodes find out if, in the next solve step, they need to solve a subproblem different from the one permanently assigned to them. This may be the case if their neighborhood currently contains a failing node $\mathcal{C}_i$. Since the partner $\mathcal{C}^\ast_i$
of the failing node will have such a compute node in its extended neighborhood, each $\mathcal{C}^\ast_i$ first checks if the associated $\mathcal{C}_i$ is back to work. If yes, then $\mathcal{C}_i$ collects all relevant static and current dynamic data for the $V_i$ subproblem as well as redundantly stored information originally assigned to it  from the corresponding neighbors, and acts in the next cycle
as normal compute node until a new failure occurs. Also, $\mathcal{C}^\ast_i$ gives up its role as partner of a failed node, and returns to normal mode.  If no, then the partner $\mathcal{C}^\ast_i$ of the still failing node $\mathcal{C}_i$ selects uniformly at random an integer $s$ from $\{0,1,\ldots,l\}$. If $s\neq 0$, then in the next solve step the neighbor $\mathcal{C}_{j_s}$ is charged with solving the $V_i$ subproblem normally assigned to the failing node $\mathcal{C}_i$ instead of its own. No further action is needed.\\
Finally, during this communication step, information about $R_0r$ is computed in a distributed way, and communicated together with the scalars $e_i$
from Step 1 to $\mathcal{S}_0$. 
\item[5] {\it Continuation/Termination}. The usual criteria for termination are based on the error indicator $\epsilon$, on the iteration count, or on the elapsed time, and are checked at the server $\mathcal{S}_0$. If none of them is satisfied, then return to Step 1. 
\ei

\medskip
{\bf Remarks}:
\bi
\item
Concerning the error indicator 
$$
\epsilon^2 =\sum_{i=0}^n e_i,
$$
we clarify its meaning as follows. Roughly speaking, $\epsilon^2$, when evaluated in Step 1 of a cycle corresponding to the recursion step 
(\ref{Rec0}), is not the value of 
$$
\epsilon_{m+1}^2:=\bar{\lambda}^{-1}a(Pe^{(m+1)},e^{(m+1)})=\sum_{i=0}^n \xi_i a_i(T_ie^{(m+1)},T_ie^{(m+1)}),
$$
which according to (\ref{NE3}) relates to the squared energy norm error of $u^{(m+1)}$ since
$$
(\lambda_{\min}/\bar{\lambda})\|e^{(m+1)}\|^2 \le \epsilon_{m+1}^2 \le (\lambda_{\max}/\bar{\lambda})\|e^{(m+1)}\|^2).
$$
However, we have 
$$
\epsilon^2= \xi_0 a_0(T_0e^{(m+1)},T_0e^{(m+1)}) +\sum_{i=1}^n \xi_i a_i(T_ie^{(m)},T_ie^{(m)}),
$$
which is almost identical with $\epsilon_{m}^2$. Thus,  $\epsilon$ is a good error measure for the previous iterate $u^{(m)}$ if $\kappa\approx \bar{\kappa}$ is moderate. Indeed, in  matrix-vector notation we have 
$$
\xi_i a_i(T_ie^{(m)},T_ie^{(m)}) = \xi_i  (A_i^{-1}R_i^Tr)^T A_i  (A_i^{-1}R_i^Tr)=\xi_i r_i^T A_i^{-1}r_i=r_i^T d_i =e'_i
$$
for any $i=0,1,\ldots,n$ (recall that $\epsilon$ is computed in Step 1 at the special server $\mathcal{S}_0$ using the value $e_0$ from the current cycle but the values $e_i$,
$i=1,\ldots,n$, from the previous cycle).

\item The numerical experiments for the one-step method (\ref{Rec0}) reported in subsections \ref{sec13} and \ref{sec14} used the steepest descent rule (\ref{XiOpt}). The computation of the required  additional global scalar products 
needs to be performed using additional communication with $\mathcal{S}_0$ during Step 2. 

\item When implementing the accelerated Schwarz iteration (\ref{u}-\ref{v}), only slight changes are required. For sure, we have to store, maintain, and communicate distributed data arrays for two vectors $x$ and $y$ representing the $u$ and $v$ iterates, respectively, and for the residual vector associated with their linear combination $z$ representing the $w$ iterate,
compare (\ref{u}). The corresponding computations can be subsumed in Steps 1 and 3, the data size in the communication steps is increased by a factor $3/2$. 

A new challenge is to provide 
an estimate for the number of active compute nodes $p_m$ which enters the parameters $\alpha_m,\beta_m$ needed in Step 1 and 3. Unfortunately,
the value of $p_m$ becomes available only after Step 2, as it may have changed due to newly failing and recovered compute nodes. Since the new fails are detected locally, counting them can only be done at $\mathcal{S}_0$. This can be achieved during Step 2
by first sending a bit from each newly created partner node through its connection to $\mathcal{S}_0$ where they are counted. 
The overall count of new faults is returned to all active nodes together with the data for the coarse problem update and $\epsilon$, and subtracted from the previous count of active nodes available at each node. This value is again corrected by adding the number of recovered compute nodes in Step 4
(to achieve their count requires communication with the specialized compute node as discussed for Step 2). In other words, the correct value
of $p_m$ is available only after the bulk of the computation for the recursion step (\ref{Rec0}) is performed in Step 1 and 3 of the associated cycle.
As mentioned in section \ref{sec1}, there are several remedies such as performing linear combinations needed for the $u$ and $v$ updates
in an extra update step after Step 4 which requires extra work and storage, especially in Step 1, by maintaining distributed vector representations separately for the residuals of the $u$ and $v$ iterates.
A cheaper alternative is to determine the parameters $\alpha_m,\beta_m$ from (\ref{AM}), (\ref{BM}) with $p_m$ replaced by the already available value for $p_{m-1}$, or a lower bound
for the latter. According to our preliminary numerical experiments for the accelerated method, having a correct value for $p_m$ may matter only if a significant part of the network is failing.

\ei

\end{document}